\documentclass[10pt]{article}
\usepackage{amssymb,amsmath}
\usepackage{latexsym}
\numberwithin{equation}{section}

\newtheorem{theorem+}           {Theorem}      [section]
\newtheorem{definition+} {Definition}      [section]
\newtheorem{lemma+}  {Lemma}  [section]
\newtheorem{corollary+}  {Corollary} [section]
\newtheorem{proposition+}  {Proposition} [section]
\newtheorem{example+}  {Example}  [section]
\newtheorem{remark+}  {Remark}  [section]
\newtheorem{problem+}  {Problem}  [section]

\newenvironment{theorem}{\begin{theorem+}\sl}{\end{theorem+}\rm}
\newenvironment{problem}{\begin{problem+}\sl}{\end{problem+}\rm}
\newenvironment{definition}{\begin{definition+}\rm}{\end{definition+}\rm}
\newenvironment{lemma}{\begin{lemma+}\sl}{\end{lemma+}\rm}
\newenvironment{corollary}{\begin{corollary+}\sl}{\end{corollary+}\rm}

\newenvironment{remark}{\begin{remark+}\rm}{\end{remark+}\rm}
\newenvironment{proof}{\medbreak\noindent{\it Proof.\ }\rm}{\hfill$\square$\medbreak\rm}

\usepackage{mystyle}
\begin{document}

\begin{center}
{\Large\bf Equilibrium problems for infinite dimensional\\ vector
potentials  with external fields}
\end{center}

\begin{center}
{\large Natalia Zorii}
\end{center}

\begin{abstract} The study deals with a minimal energy problem in the presence of an
external field~$\mathbf{f}=(f_i)_{i\in I}$ over noncompact classes
of vector measures $\mu=(\mu^i)_{i\in I}$ of infinite dimension in a
locally compact space. The components $\mu^i$ are positive measures
(charges) normalized by $\int g_i\,d\mu^i=a_i$ (where $a_i$ and
$g_i$ are given) and supported by given closed sets~$A_i$ with the
sign~$+1$ or~$-1$ prescribed such that $A_i\cap A_j=\varnothing$
whenever ${\rm sign}\,A_i\ne{\rm sign}\,A_j$, and the law of
interaction of $\mu^i$, $i\in I$, is determined by the interaction
matrix $\bigl({\rm sign}\,A_i\,{\rm sign}\,A_j\bigr)_{i,j\in I}$.
For all positive definite kernels satisfying Fuglede's condition of
consistency between the vague ($={}$weak$*$) and strong topologies,
sufficient conditions for the existence of equilibrium measures are
established and properties of their uniqueness, vague compactness,
and continuity under exhaustion of~$A_i$ by compact~$K_i$ are
studied. We also obtain variational inequalities for the
$\mathbf{f}$-weighted equilibrium potentials, single out their
characteristic properties, and analyze continuity of the equilibrium
constants.


\end{abstract}

\section{Introduction}
\label{intro} The interest to minimal energy problems in the
presence of an external field, initially inspired by
C.\,F.~Gauss~\cite{Gauss} and further experiencing a new growth due
to work of~O.~Frostman~\cite{Fr} and Polish and Japanese
mathematicians (F.~Leja, J.~G\'{o}rski, W.~Kleiner, J.~Siciak and
S.~Kametani, M.~Ohtsuka, N.~Ninomiya; see~\cite{Leja,O} and the
references cited therein), has been motivated by their direct
relations with the Dirichlet and balayage problems.

A new impulse to this part of potential theory (which is often
referred to as the Gauss variational problem) came in the 1980's
when A.\,A.~Gonchar and E.\,A.~Rakhmanov~\cite{GR0,GR1},
H.\,N.~Mhaskar and E.\,B.~Saff~\cite{MaS} efficiently applied
logarithmic potentials with external fields in the investigation of
orthogonal polynomials and rational approximations to analytic
functions; for references to subsequent publications, see the
books~\cite{NS,ST}.

We shall consider the Gauss variational problem in a rather general
setting, over classes of vector measures of infinite dimension in a
locally compact Hausdorff space~$\mathrm X$. In case the measures
are of finite dimension, the vector setting of the problem goes back
to~\cite[\S~2.9]{O}; see also~\cite{GR0,GR}, related to the
logarithmic kernel in the plane. To formulate the problem and
shortly outline the results obtained, we start by introducing
briefly relevant notions.

Let $\mathfrak M=\mathfrak M(\mathrm X)$ denote the linear space of
all real-valued scalar Radon measures~$\nu$ on~$\mathrm X$ equipped
with the {\it vague\/} ($={}${\it weak}$*$) topology, i.\,e., the
topology of pointwise convergence on the class $\mathrm C_0(\mathrm
X)$ of all real-valued continuous functions~$\varphi$ on~$\mathrm X$
with compact support. A {\it kernel\/}~$\kappa$ on $\mathrm X$ is
meant to be an element from $\mathrm\Phi(\mathrm X\times\mathrm X)$,
where $\mathrm\Phi(\mathrm Y)$ consist of all lower semicontinuous
functions~$\psi:\mathrm Y\to(-\infty,\infty]$ such that
$\psi\geqslant0$ unless $\mathrm Y$ is compact.

Given $\nu,\,\nu_1\in\mathfrak M$, the {\it mutual energy\/} and the
{\it potential\/} with respect to a kernel~$\kappa$ are defined
respectively by
$$\kappa(\nu,\nu_1):=\int\kappa(x,y)\,d(\nu\otimes\nu_1)(x,y)\quad\mbox{and}\quad
\kappa(\,\cdot\,,\nu):=\int\kappa(\,\cdot\,,y)\,d\nu(y).$$ (Here and
in the sequel, when introducing notation, we shall always tacitly
assume the corresponding object on the right to be well defined.)
For $\nu=\nu_1$ the mutual energy $\kappa(\nu,\nu_1)$ gives the {\it
energy\/} of~$\nu$. The set of all $\nu\in\mathfrak M$ with
$-\infty<\kappa(\nu,\nu)<\infty$ will be denoted by~$\mathcal
E=\mathcal E_\kappa$.

We shall be mainly concerned with a {\it positive definite\/}
kernel~$\kappa$, which means that it is symmetric (i.\,e.,
$\kappa(x,y)=\kappa(y,x)$ for all $x,\,y\in\mathrm X$) and the
energy $\kappa(\nu,\nu)$, $\nu\in\mathfrak M$, is nonnegative
whenever defined. Then $\mathcal E$ forms a pre-Hil\-bert space with
the scalar product $\kappa(\nu,\nu_1)$ and the seminorm
$\|\nu\|_\mathcal E:=\sqrt{\kappa(\nu,\nu)}$ (see~\cite{F1}). A
positive definite kernel~$\kappa$ is called {\it strictly positive
definite\/} if the seminorm $\|\,\cdot\,\|_\mathcal E$ is a norm.

Given a closed set $E\subset\mathrm X$, let $\mathfrak M^+(E)$
consist of all nonnegative measures $\nu\in\mathfrak M$ supported
by~$E$, and let $\mathcal E^+(E):=\mathfrak M^+(E)\cap\mathcal E$.
Also write $\mathfrak M^+:=\mathfrak M^+(\mathrm X)$ and $\mathcal
E^+:=\mathcal E^+(\mathrm X)$.

We consider a countable, locally finite collection
$\mathbf{A}=(A_i)_{i\in I}$ of fixed closed sets $A_i\subset\mathrm
X$ with the sign~$+1$ or $-1$ prescribed such that the oppositely
signed sets are mutually disjoint. Let $\mathfrak M(\mathbf{A})$
stand for the Cartesian product $\prod_{i\in I}\,\mathfrak
M^+(A_i)$; then an element~$\mu$ of~$\mathfrak M(\mathbf{A})$ is a
vector measure $(\mu^i)_{i\in I}$ with the components
$\mu^i\in\mathfrak M^+(A_i)$. If, moreover, $\mathbf{u}=(u_i)_{i\in
I}$ is a vector-valued function, we shall write
$\langle\mathbf{u},\mu\rangle:=\sum_{i\in I}\,\int u_i\,d\mu^i$.

Let a kernel~$\kappa$ be fixed. Corresponding to an electrostatic
interpretation, we assume that the interaction of point charges
lying on the conductors~$A_i$, $i\in I$, is characterized by the
interaction matrix $(\alpha_i\alpha_j)_{i,j\in I}$, where
$\alpha_i:={\rm sign}\,A_i$. Given vector measures
$\mu,\,\mu_1\in\mathfrak M(\mathbf{A})$, we define the {\it mutual
energy\/}
\begin{equation}\label{vectoren}\kappa(\mu,\mu_1):=\sum_{i,j\in
I}\,\alpha_i\alpha_j\kappa(\mu^i,\mu_1^j)\end{equation} and the {\it
vector potential\/} $\kappa_\mu(x)$, $x\in\mathrm X$, as a
vector-valued function with the components
\begin{equation}\label{vectorpot}\kappa^i_\mu(x):=\sum_{j\in
I}\,\alpha_i\alpha_j\kappa(x,\mu^j),\quad i\in I.\end{equation} For
$\mu=\mu_1$ the mutual energy $\kappa(\mu,\mu_1)$ defines the {\it
energy\/} of~$\mu$. Let $\mathcal E(\mathbf{A})$ consist of all
$\mu\in\mathfrak M(\mathbf{A})$ whose energy $\kappa(\mu,\mu)$ is
finite.

Fix also a vector-valued function $\mathbf{f}=(f_i)_{i\in I}$ to be
treated as an external field. The $\mathbf{f}$-{\it weight\-ed
vector potential\/} and the $\mathbf{f}$-{\it weighted energy\/} of
$\mu\in\mathcal E(\mathbf{A})$ are then defined by
\begin{align}\label{wpot}\mathbf{W}_\mu&:=\kappa_\mu+\mathbf{f},\\
\label{wen}G_{\mathbf{f}}(\mu)&:=
\kappa(\mu,\mu)+2\langle\mathbf{f},\mu\rangle,\end{align}
respectively. In the present study we shall be mainly focused with
the case where either $f_i\in\mathrm\Phi(\mathrm X)$ for all $i\in
I$, or $f_i=\alpha_i\kappa(\,\cdot\,,\sigma)$, $i\in I$ (here
$\sigma\in\mathcal E$ is given).

We also fix a numerical vector $\mathbf{a}=(a_i)_{i\in I}$ with
$a_i>0$ for all~$i\in I$ and a vector-valued function
$\mathbf{g}=(g_i)_{i\in I}$, where $g_i:A_i\to(0,\infty)$ are
continuous. We shall be interested in the problem of minimizing
$G_{\mathbf{f}}(\mu)$ over the class of all $\mu\in\mathcal
E(\mathbf{A})$ with $\langle g_i,\mu^i\rangle=a_i$, $i\in I$.

The main question is whether equilibrium
measures~$\lambda_{\mathbf{A}}$ in the minimal $\mathbf{f}$-weighted
energy problem exist. If $\mathbf{A}$ is finite, $A_i$ is compact
and $f_i\in\mathrm\Phi(\mathrm X)$ for every $i\in I$, while
$\kappa(x,y)$ is continuous on $A_i\times A_j$ whenever
$\alpha_i\ne\alpha_j$, then the existence of
those~$\lambda_{\mathbf{A}}$ can easily be established by exploiting
the vague topology only (see~\cite{O};
cf.~also~\cite{GR0,GR,NS,ST}). However, the question becomes rather
nontrivial if any of these four assumptions is dropped.

To solve the problem on the existence of equilibrium
measures~$\lambda_{\mathbf{A}}$ in the general case where
$\mathbf{A}$ is infinite and (or) $A_i$, $i\in I$, are noncompact,
we restrict ourselves to positive definite kernels~$\kappa$ and work
out an approach based on the following arguments.

The set $\mathcal E(\mathbf{A})$ is shown to be a semimetric space
with the semimetric (see Sect.~\ref{sec:semimetric})
\begin{equation}\label{vseminorm}
\|\mu_1-\mu_2\|_{\mathcal E(\mathbf{A})}:=\Bigl[\sum_{i,j\in
I}\,\alpha_i\alpha_j\kappa(\mu^i_1-\mu^i_2,\mu^j_1-\mu^j_2)\Bigr]^{1/2},
\end{equation}
and one can define an inclusion $R$ of $\mathcal E(\mathbf{A})$ into
the pre-Hilbert space~$\mathcal E$ such that $\mathcal
E(\mathbf{A})$ is isometric to its $R$-image, the latter being
regarded as a semimetric subspace of~$\mathcal E$.

Another crucial fact is that, for rather general $\kappa$,
$\mathbf{g}$, and~$\mathbf{a}$, the topological subspace of
$\mathcal E(\mathbf{A})$ consisting of all $\mu$ with $\langle
g_i,\mu^i\rangle\leqslant a_i$, $i\in I$, turns out to be complete
(see Theorem~\ref{th:strong}).

Using these arguments, we obtain sufficient conditions for the
existence of equilibrium measures~$\lambda_{\mathbf{A}}$ and
establish statements on their uniqueness and vague compactness (see
Lemma~\ref{lemma:unique} and Theorem~\ref{exist}). Continuity
properties of equilibrium measures under exhaustion of~$\mathbf{A}$
by~$\mathbf{K}$ with compact~$K_i$, $i\in I$, are analyzed as well
(see Theorem~\ref{cor:cont}).

We also establish variational inequalities for the
$\mathbf{f}$-weighted equilibrium potentials
$\mathbf{W}_{\lambda_{\mathbf{A}}}$ (see Theorems~\ref{th:descpot1}
and~\ref{th:descpot2}); some of those inequalities are shown to be
characteristic (see~Theorem~\ref{th:char}). In particular, there
exist numbers $C^i_{\mathbf{A}}$, $i\in I$, called the
$\mathbf{f}$-{\it weighted equilibrium constants\/}, such that
\begin{align*}
a_i\,W_{\lambda_{\mathbf{A}}}^i(x)&\geqslant C^i_{\mathbf{A}}\,g(x)\quad\mbox{n.\,e. in \ } A_i,\\
G_{\mathbf{f}}(\lambda_{\mathbf{A}})&\leqslant\sum_{i\in
I}\,C^i_{\mathbf{A}}+\langle\mathbf{f},\lambda_{\mathbf{A}}\rangle,
\end{align*}
where {\it n.\,e.\/}~({\it nearly everywhere\/}) means that the set
of all $x\in A_i$ for which the inequality fails to hold has
interior capacity zero; and these inequalities determine uniquely
equilibrium measures among all the admissible ones. Under proper
additional restrictions, it is also true that
\[a_i\,W_{\lambda_{\mathbf{A}}}^i(x)\leqslant
C^i_{\mathbf{A}}\,g(x)\quad\mbox{for all \ }x\in
S(\lambda_{\mathbf{A}}^i).\] The equilibrium constants are uniquely
determined and can be written in either of the forms
\[C^i_{\mathbf{A}}=\bigl\langle
W^i_{\lambda_{\mathbf{A}}},\lambda_{\mathbf{A}}^i\bigr\rangle="\!\inf_{x\in
A_i}\!"\,\,\frac{a_i\,W_{\lambda_{\mathbf{A}}}^i(x)}{g(x)},\] the
infimum being taken over all~$A_i$ excepting probably its subset of
interior capacity zero. Furthermore, for rather general $\kappa$,
$\mathbf{g}$, $\mathbf{a}$, and~$\mathbf{f}$, these constants are
shown to be continuous under exhaustion of~$\mathbf{A}$
by~$\mathbf{K}$ with compact~$K_i$, $i\in I$
(see~Theorem~\ref{cor:cont}).

The results obtained and the approach applied develop and generalize
the corresponding ones from the author's
articles~\cite{Z4,Z5a,Z5,Z6}, related to vector measures of finite
dimensions.

\section{Preliminaries: topologies, consistent and perfect
kernels}\label{sec:2}

In all that follows, we shall always suppose the kernel~$\kappa$ to
be positive definite. In addition to the {\it strong\/} topology
on~$\mathcal E$, determined by the seminorm
$\|\nu\|:=\|\nu\|_\mathcal E$, it is often useful to consider the
{\it weak\/} topology on~$\mathcal E$, defined by means of the
seminorms $\nu\mapsto|\kappa(\nu,\mu)|$, $\mu\in\mathcal E$
(see~\cite{F1}). The Cauchy-Schwarz inequality
\begin{equation*}
|\kappa(\nu,\mu)|\leqslant\|\nu\|\,\|\mu\|,\quad\mbox{where
$\nu,\,\mu\in\mathcal E$,}\end{equation*} implies immediately that
the strong topology on $\mathcal E$ is finer than the weak one.

In~\cite{F1,F2}, B.~Fuglede introduced the following two {\it
equivalent\/} properties of consistency between the induced strong,
weak, and vague topologies on~$\mathcal E^+$:
\begin{itemize}
\item[\rm(C$_1$)] {\it Every strong Cauchy net in
$\mathcal E^+$ converges strongly to every its vague cluster
point;}\smallskip
\item[\rm(C$_2$)] {\it Every strongly bounded and vaguely convergent net in
$\mathcal E^+$ converges weakly to the vague limit.}
\end{itemize}

\begin{definition}
Following Fuglede~\cite{F1}, we call a kernel~$\kappa$ {\it
consistent} if it satisfies either of the properties~(C$_1$)
and~(C$_2$), and {\it perfect\/} if, in addition, it is strictly
positive definite.\end{definition}

\begin{remark} One has to consider {\it nets\/} or {\it filters\/}
in~$\mathfrak M^+$ instead of sequences, since the vague topology in
general does not satisfy the first axiom of countability. We follow
Moore's and Smith's theory of convergence, based on the concept of
nets (see~\cite{MS}; cf.~also~\cite[Chap.~0]{E2} and
\cite[Chap.~2]{K}). However, if $\mathrm X$ is metrizable and
countable at infinity, then $\mathfrak M^+$ satisfies the first
axiom of countability (see~\cite[Lemma~1.2.1]{F1}) and the use of
nets may be avoided.\end{remark}

\begin{theorem}{\rm(Fuglede \cite{F1})}\label{th:1} A kernel $\kappa$ is perfect if and
only if $\mathcal E^+$ is strongly complete and the strong topology
on~$\mathcal E^+$ is finer than the vague one.\end{theorem}

\begin{remark} In $\mathbb R^n$, $n\geqslant 3$, the
Newtonian kernel $|x-y|^{2-n}$ is perfect~\cite{Car}. So are the
Riesz kernel $|x-y|^{\alpha-n}$, $0<\alpha<n$, in~$\mathbb R^n$,
$n\geqslant2$~\cite{D1,D2}, and the restriction of the kernel
$-\log\,|x-y|$ in~$\mathbb R^2$ to an open unit ball~\cite{L}.
Furthermore, if $D$ is an open set in~$\mathbb R^n$, $n\geqslant 2$,
and its generalized Green function~$g_D$ exists (see,
e.\,g.,~\cite[Th.~5.24]{HK}), then $g_D$ is perfect as
well~\cite{E1}.\end{remark}

\begin{remark} As is seen from the above definitions and Theorem~\ref{th:1}, the concept of consistent
or perfect kernels is an efficient tool in minimal energy problems
over classes of {\it nonnegative scalar\/} Radon measures with
finite energy. Indeed, the theory of capacities of {\it sets\/} has
been developed in~\cite{F1} exactly for those kernels. We shall show
below that this concept is efficient, as well, in minimal energy
problems over classes of {\it vector measures\/} of finite or
infinite dimensions. This is guaranteed by a theorem on the
completeness of proper subspaces of the semimetric space~$\mathcal
E(\mathbf{A})$, to be stated in Sect.~\ref{sec:strong}.\end{remark}

\section{Condensers. Vector measures; their energies and potentials}

\subsection{Condensers of countably many plates. Associated vector measures}
Let $I^+$ and $I^-$ be countable (finite or infinite) disjoint sets
of indices~$i\in\mathbb N$, where the latter is allowed to be empty,
and let $I$ denote their union. Assume that to every $i\in I$ there
corresponds a nonempty, closed set~$A_i\subset\mathrm X$.

\begin{definition} A collection $\mathbf{A}=(A_i)_{i\in I}$ is
called an $(I^+,I^-)$-{\it condenser\/} (or simply a {\it
condenser\/}) in~$\mathrm X$ if every compact subset of~$\mathrm X$
intersects with at most finitely many~$A_i$ and
\begin{equation}
A_i\cap A_j=\varnothing\quad\mbox{for all \ } i\in I^+, \ j\in I^-.
\label{non}
\end{equation}\end{definition}

The sets $A_i$, $i\in I^+$, and $A_j$, $j\in I^-$, are called the
{\it positive\/} and, respectively, {\it negative plates\/} of the
condenser~$\mathbf{A}$. Note that any two equally sign\-ed plates
can intersect each other.

Given $I^+$ and $I^-$, let $\mathfrak C=\mathfrak C(I^+,I^-)$ be the
class of all $(I^+,I^-)$-condensers in~$\mathrm X$. A condenser
$\mathbf{A}\in\mathfrak C$ will be called {\it compact\/} if so are
all~$A_i$, $i\in I$, and {\it finite\/} if $I$ is finite.

In the sequel, also the following notation will be used:
\begin{equation*}
A^+:=\bigcup_{i\in I^+}\,A_i,\qquad A^-:=\bigcup_{i\in I^-}\,A_i.
\end{equation*} Observe that $A^+$
and~$A^-$ might both be noncompact even for a compact~$\mathbf{A}$.

Given $\mathbf{A}\in\mathfrak C$, let $\mathfrak M(\mathbf{A})$
consist of all {\it vector measures\/} $\mu=(\mu^i)_{i\in I}$, where
$\mu^i\in\mathfrak M^+(A_i)$ for all $i\in I$; that is, $\mathfrak
M(\mathbf{A})$ stands for the Cartesian product $\prod_{i\in
I}\,\mathfrak M^+(A_i)$. The product topology on~$\mathfrak
M(\mathbf{A})$, where every $\mathfrak M^+(A_i)$ is equipped with
the vague topology, will be called the $\mathbf{A}$-{\it vague
topology\/}. Since $\mathfrak M(\mathrm X)$ is Hausdorff, so is
$\mathfrak M(\mathbf{A})$ (cf.~\cite[Chap.~3, Th.~5]{K}).

A set $\mathfrak F\subset\mathfrak M(\mathbf{A})$ is called
$\mathbf{A}$-{\it vaguely bounded\/} if, for all $\varphi\in\mathrm
C_0(\mathrm X)$ and $i\in I$,
\[\sup_{\mu\in\mathfrak F}\,|\mu^i(\varphi)|<\infty.\]

\begin{lemma}\label{lem:vaguecomp} If $\mathfrak F\subset\mathfrak
M(\mathbf{A})$ is $\mathbf{A}$-vaguely bounded, then it is
$\mathbf{A}$-vaguely relatively compact.\end{lemma}

\begin{proof} Since by~\cite[Chap.~III, \S~2, Prop.~9]{B2} any
vaguely bounded part of~$\mathfrak M$ is vaguely relatively compact,
the lemma follows immediately from Tychonoff's theorem on the
product of compact spaces (see, e.\,g.,~\cite[Chap.~5,
Th.~13]{K}).
\end{proof}

\subsection{Mapping $R:\mathfrak M(\mathbf{A})\to\mathfrak M$.
Relation of $R$-equivalency on $\mathfrak
M(\mathbf{A})$}\label{sect:R}

Since each compact subset of~$\mathrm X$ intersects with at most
finitely many~$A_i$, for every $\varphi\in\mathrm C_0(\mathrm X)$
only a finite number of~$\mu^i(\varphi)$ (where $\mu\in\mathfrak
M(\mathbf{A})$ is given) are nonzero. This yields that to every
vector measure $\mu\in\mathfrak M(\mathbf{A})$ there corresponds a
unique scalar Radon measure~$R\mu\in\mathfrak M$ such that
\[R\mu(\varphi)=\sum_{i\in I}\,\alpha_i\mu^i(\varphi)\quad\mbox{for
all \ }\varphi\in\mathrm C_0(\mathrm X);\] because of~(\ref{non}),
positive and negative parts in Jordan's decomposition of $R\mu$ can
respectively be written in the form
\[ R\mu^+=\sum_{i\in I^+}\,\mu^i,\qquad R\mu^-=\sum_{i\in
I^-}\,\mu^i.\]

Of course, the inclusion $\mathfrak M(\mathbf{A})\to\mathfrak M$
thus defined is in general non-injective, i.\,e., one may choose
$\mu_1,\,\mu_2\in\mathfrak M(\mathbf{A})$ so that $\mu_1\ne\mu_2$,
while $R\mu_1=R\mu_2$. We shall call $\mu_1,\,\mu_2\in\mathfrak
M(\mathbf{A})$ {\it $R$-eq\-uivalent\/} if $R\mu_1=R\mu_2$ --- or,
which is equivalent, whenever $\sum_{i\in I}\,\mu_1^i=\sum_{i\in
I}\,\mu_2^i$.

Observe that the relation of $R$-equivalency implies that of
identity (and, hence, these two relations on~$\mathfrak
M(\mathbf{A})$ are actually equivalent) if and only if all~$A_i$,
$i\in I$, are mutually disjoint.

\begin{lemma}\label{lem:vague'}The $\mathbf{A}$-vague convergence of
$(\mu_s)_{s\in S}\subset\mathfrak M(\mathbf{A})$
to~$\mu_0\in\mathfrak M(\mathbf{A})$ implies the vague convergence
of $(R\mu_s)_{s\in S}$ to~$R\mu_0$.\end{lemma}

\begin{proof} This is obvious in view of the fact that the
support of any $\varphi\in\mathrm C_0(\mathrm X)$ might have points
in common with only finitely many~$A_i$.
\end{proof}

\begin{remark} Lemma~\ref{lem:vague'} in general can not be
inverted. However, if all $A_i$, $i\in I$, are mutually disjoint,
then the vague convergence of $(R\mu_s)_{s\in S}$ to~$R\mu_0$
implies the $\mathbf{A}$-vague convergence of $(\mu_s)_{s\in S}$
to~$\mu_0$. This can be seen by using the Tietze-Urysohn extension
theorem.\end{remark}

\subsection{Energies and potentials of vector measures and their $R$-images}
In accordance with an electrostatic interpretation of a
condenser~$\mathbf{A}$, we suppose that the law of interaction of
charges lying on its plates~$A_i$, $i\in I$, is determined by the
interaction matrix $(\alpha_i\alpha_j)_{i,j\in I}$, where
\[\alpha_i:=\left\{
\begin{array}{rll} +1 & \mbox{if} & i\in I^+,\\ -1 & \mbox{if} & i\in
I^-.\\ \end{array} \right.\] Given vector measures
$\mu,\,\mu_1\in\mathfrak M(\mathbf{A})$, we define the {\it mutual
energy\/} $\kappa(\mu,\mu_1)$ and the {\it vector potential\/}
$\kappa_\mu=(\kappa^i_\mu)_{i\in I}$ by~(\ref{vectoren})
and~(\ref{vectorpot}), respectively. If $\mu=\mu_1$, then
$\kappa(\mu,\mu_1)$ defines the {\it energy\/} $\kappa(\mu,\mu)$ of
$\mu$.

\begin{lemma}\label{enfinite}  For $\mu\in\mathfrak M(\mathbf{A})$ to
be of finite energy, it is necessary and sufficient that
$\mu^i\in\mathcal E$ for all $i\in I$ and \[\sum_{i\in
I}\,\|\mu^i\|^2<\infty.\]\end{lemma}

\begin{proof} This follows  immediately from the definition of
$\kappa(\mu,\mu)$ in view of the inequality
$2\kappa(\nu_1,\nu_2)\leqslant\|\nu_1\|^2+\|\nu_2\|^2$ for
$\nu_1,\,\nu_2\in\mathcal E$.
\end{proof}

To establish relations between energies and potentials of vector
measures $\mu\in\mathfrak M(\mathbf{A})$ and those of their (scalar)
$R$-images $R\mu\in\mathfrak M$, we start with the following two
lemmas, the first one being well known (see, e.\,g., \cite{F1}).

\begin{lemma}\label{lemma:lower}
If $\mathrm Y$ is a locally compact Hausdorff space and
$\psi\in\mathrm\Phi(\mathrm Y)$ is given, then the map
$\nu\mapsto\langle\psi,\nu\rangle$ is vaguely lower semicontinuous
on $\mathfrak M^+(\mathrm Y)$.\end{lemma}

\begin{lemma}\label{integral} Fix $\mu\in\mathfrak M(\mathbf{A})$ and
$\psi\in\mathrm\Phi(\mathrm X)$. If $\langle \psi,R\mu\rangle$ is
well defined, then
\begin{equation}
\langle\psi,R\mu\rangle=\sum_{i\in I}\,\alpha_i\langle\psi,
\mu^i\rangle,\label{lemma11}
\end{equation}
and $\langle\psi,R\mu\rangle$ is finite if and only if the series on
the right converges absolutely.\end{lemma}

\begin{proof} We can assume $\psi$ to be nonnegative,
for if not, we replace $\psi$ by a function~$\psi'\geqslant0$
obtained by adding to~$\psi$ a suitable constant~$c>0$, which is
always possible since a lower semicontinuous function is bounded
from below on a compact space. Hence,
\[\langle\psi,R\mu^+\rangle\geqslant\sum_{i\in I^+, \ i\leqslant N}\,\langle
\psi,\mu^i\rangle\quad\mbox{for all \ }N\in\mathbb N.\] On the other
hand, the sum of $\mu^i$ over all $i\in I^+$ that do not exceed~$N$
approaches~$R\mu^+$ vaguely as $N\to\infty$; consequently, by
Lemma~\ref{lemma:lower},
\begin{equation*}
\langle\psi,R\mu^+\rangle\leqslant\lim_{N\to\infty}\,\sum_{i\in I^+,
\ i\leqslant N}\,\langle\psi,\mu^i\rangle.\end{equation*} Combining
the last two inequalities and then letting $N\to\infty$ yields
\[\langle\psi,R\mu^+\rangle=\sum_{i\in I^+}\,\langle\psi,\mu^i\rangle.\] Since the same
holds true for $R\mu^-$ and~$I^-$ instead of~$R\mu^+$ and~$I^+$, the
lemma follows.
\end{proof}

\begin{corollary}\label{pot.ener} Fix $\mu,\,\mu_1\in\mathfrak
M(\mathbf{A})$ and $x\in\mathbf X$. Then
\begin{align}
\kappa(R\mu,R\mu_1)&=\sum_{i,j\in
I}\,\alpha_i\alpha_j\kappa(\mu^i,\mu_1^j),\label{mutual}\\
\kappa(x,R\mu)&=\sum_{i\in I}\,\alpha_i\kappa(x,\mu^i),\label{poten}
\end{align}
each of the identities being understood in the sense that either of
its sides is well defined whenever so is the other one and then they
coincide. Furthermore, the left-hand side in~{\rm(\ref{mutual})} or
in~{\rm(\ref{poten})} is finite if and only if the corresponding
series on the right converges absolutely.\end{corollary}

\begin{proof} Relation (\ref{poten}) is a direct consequence
of~(\ref{lemma11}), while (\ref{mutual}) follows from Fubini's
theorem (cf.~\cite[\S~8, Th.~1]{B3}) and Lemma~\ref{integral} on
account of the fact that $\kappa(x,\nu)$, where $\nu\in\mathfrak
M^+$ is given, is lower semicontinuous on~$\mathrm X$ (see,
e.\,g.,~\cite{F1}).
\end{proof}

When comparing (\ref{vectoren}) and (\ref{vectorpot}) with
(\ref{mutual}) and (\ref{poten}), respectively, we obtain

\begin{corollary}\label{relation}Given $\mu,\,\mu_1\in\mathfrak
M(\mathbf{A})$, $x\in\mathrm X$, and $i\in I$,
\begin{align}
\kappa(\mu,\mu_1)&=\kappa(R\mu,R\mu_1),\label{Ren}\\
\kappa^i_\mu(x)&=\alpha_i\kappa(x,R\mu).\label{Rpot}
\end{align}
\end{corollary}

\subsection{Semimetric space of vector measures of finite
energy}\label{sec:semimetric}

Let $\mathcal E(\mathbf{A})$ consist of all $\mu\in\mathfrak
M(\mathbf{A})$ with finite energy $\kappa(\mu,\mu)$. Since
$\mathfrak M(\mathbf{A})$ is a convex cone, it follows from
Lemma~\ref{enfinite} that so is $\mathcal E(\mathbf{A})$.

\begin{lemma}\label{lemma:semimetric}The cone $\mathcal E(\mathbf{A})$ forms a semimetric space with
the semimetric $\|\,\cdot\,\|_{\mathcal E(\mathbf{A})}$ defined
by~{\rm(\ref{vseminorm})}. This semimetric is a metric if and only
if the kernel $\kappa$ is strictly positive definite while all
$A_i$, $i\in I$, are mutually disjoint.
\end{lemma}

\begin{proof}Fix $\mu_1,\,\mu_2\in\mathcal E(\mathbf{A})$. Applying Corollary~\ref{pot.ener} to
$\kappa(R\mu_k,R\mu_\ell)$, $k,\,\ell=1,\,2$, we get
\begin{equation*}\label{isometric}\|R\mu_1-R\mu_2\|_{\mathcal
E}^2=\sum_{i,j\in
I}\,\alpha_i\alpha_j\kappa(\mu^i_1-\mu^i_2,\mu^j_1-\mu^j_2).\end{equation*}
When compared with~(\ref{vseminorm}), this yields
\begin{equation}\label{seminorm}\|\mu_1-\mu_2\|^2_{\mathcal E(\mathbf{A})}=\|R\mu_1-R\mu_2\|_{\mathcal
E}^2.\end{equation} Since $\|\cdot\|_{\mathcal E}$ is a seminorm
on~$\mathcal E$, the proof is complete.
\end{proof}

In all that follows, $\mathcal E(\mathbf{A})$ will always be treated
as a semimetric space with the
semimetric~\mbox{$\|\,\cdot\,\|:=\|\,\cdot\,\|_{\mathcal
E(\mathbf{A})}$}. Then, by~(\ref{seminorm}), $\mathcal
E(\mathbf{A})$ and its $R$-image become isometric. Similarly with
the terminology in~$\mathcal E$, the topology on $\mathcal
E(\mathbf{A})$ will be called {\it strong}.

Two elements of~$\mathcal E(\mathbf{A})$, $\mu_1$ and~$\mu_2$, are
said to be {\it equivalent in\/}~$\mathcal E(\mathbf{A})$ if
$\|\mu_1-\mu_2\|=0$. Observe that the equivalence in~$\mathcal
E(\mathbf{A})$ implies $R$-equivalence (i.\,e., then
$R\mu_1=R\mu_2$) provided the kernel~$\kappa$ is strictly positive
definite, and it implies the identity (i.\,e., then $\mu_1=\mu_2$)
if, moreover, all~$A_i$, $i\in I$, are mutually disjoint.

A vector-valued proposition $\mathbf{u}=(u_i)_{i\in I}$ involving a
variable point $x\in\mathbf X$ is said to subsist {\it nearly
everywhere\/}~(n.\,e.) in~$E$, where $E$ is a given subset
of~$\mathbf X$, if for every $i\in I$ the set of all $x\in E$ for
which $u_i$ fails to hold is of interior capacity zero.
\begin{corollary}\label{lemma:potfinite}For every $\mu\in\mathcal E(\mathbf{A})$, $\kappa_{\mu}(x)$
is defined and finite nearly everywhere in~$\mathrm
X$.\end{corollary}
\begin{proof}This is seen from~(\ref{Ren}) and~(\ref{Rpot}) in view of the fact
that the potential $\kappa(x,\nu)$ of any $\nu\in\mathcal E$ is
defined and finite n.\,e.~in~$\mathrm X$
(see~\cite{F1}).
\end{proof}
\begin{corollary}\label{lemma:potequiv}If $\mu_1$ and~$\mu_2$ are equivalent in $\mathcal E(\mathbf{A})$,
then \[\kappa_{\mu_1}(x)=\kappa_{\mu_2}(x)\quad\mbox{n.\,e. in \ }
\mathrm X.\]
\end{corollary}
\begin{proof}Indeed, then $R\mu_1$ and~$R\mu_2$ are equivalent in $\mathcal
E$ by (\ref{seminorm}). Hence, $\kappa(x,R\mu_1)=\kappa(x,R\mu_2)$
nearly everywhere in~$\mathrm X$ (see~\cite{F1}), which together
with~(\ref{Rpot}) proves the corollary.
\end{proof}

\section{Minimal $\mathbf{f}$-weighted
energy problem}\label{sec:statement}

From now on the external field $\mathbf{f}=(f_i)_{i\in I}$ will
always be of the following structure. For every $i\in I$, there are
$f_{i1},\,f_{i2}\in\mathrm\Phi(\mathrm X)$ such that
$f_{i2}\ne\infty$ n.\,e.~in~$\mathrm X$ and
\begin{equation*}\label{difference}f_i(x)=f_{i1}(x)-f_{i2}(x),\quad x\in\mathrm X,\end{equation*}
where the value on the left is defined if and only if so is that on
the right and then they coincide. Such an $f_i$ is defined and
${}\ne-\infty$ n.\,e.~in~$\mathrm X$ and is universally measurable,
i.\,e., measurable with respect to every $\nu\in\mathfrak M$. Also
note that, for any $\mu\in\mathfrak M(\mathbf{A})$,
$\langle\mathbf{f},\mu\rangle$ is finite if and only if $\sum_{i\in
I}\,\langle f_i,\mu^i\rangle$ converges absolutely.

Given $\mu\in\mathcal E(\mathbf{A})$, we then define the
$\mathbf{f}$-{\it weighted vector potential\/} $\mathbf{W}_\mu$ and
the $\mathbf{f}$-{\it weighted energy\/} $G_{\mathbf{f}}(\mu)$
by~(\ref{wpot}) and~(\ref{wen}), respectively. Note that, according
to Corollary~\ref{lemma:potfinite}, $\mathbf{W}_\mu$ is defined and
${}\ne-\infty$ n.\,e.~in~$\mathrm X$. Also observe that,
by~(\ref{Ren}), (\ref{Rpot}), and Fubini's theorem,
\[G_{\mathbf{f}}(\mu)=\bigl\langle\mathbf{W}_\mu+\mathbf{f},\mu\bigr\rangle.\]

Having fixed also a vector-valued function $\mathbf{g}=(g_i)_{i\in
I}$, where $g_i:A_i\to(0,\infty)$, $i\in I$, are continuous, and a
numerical vector $\mathbf{a}=(a_i)_{i\in I}$ with $a_i>0$, we write
\begin{equation*}\mathfrak M(\mathbf{A},\mathbf{a},\mathbf{g}):=\bigl\{\mu\in\mathfrak
M(\mathbf{A}): \ \langle g_i,\mu^i\rangle=a_i\quad\mbox{for all \ }
i\in I\bigr\},\end{equation*}
\[\mathcal E(\mathbf{A},\mathbf{a},\mathbf{g}):=\mathfrak M(\mathbf{A},\mathbf{a},\mathbf{g})\cap\mathcal
E(\mathbf{A}),\]
\[\mathcal
E_{\mathbf{f}}(\mathbf{A},\mathbf{a},\mathbf{g}):=\bigl\{\mu\in\mathcal
E(\mathbf{A},\mathbf{a},\mathbf{g}): \ \langle
\mathbf{f},\mu\rangle\mbox{ \ is finite} \bigr\}\] and further
introduce the extremal value
\begin{equation}\label{G}G_{\mathbf{f}}(\mathbf{A},\mathbf{a},\mathbf{g}):=\inf_{\mu\in\mathcal
E_{\mathbf{f}}(\mathbf{A},\mathbf{a},\mathbf{g})}\,G_{\mathbf{f}}(\mu).
\end{equation}
In (\ref{G}), as usual, the infimum over the empty set is taken to
be~$+\infty$.

\begin{problem}If $-\infty<G_{\mathbf{f}}(\mathbf{A},\mathbf{a},\mathbf{g})<\infty$, does there exist
$\lambda=\lambda_{\mathbf{A}}\in\mathcal
E_{\mathbf{f}}(\mathbf{A},\mathbf{a},\mathbf{g})$ with
\[G_{\mathbf{f}}(\lambda)=G_{\mathbf{f}}(\mathbf{A},\mathbf{a},\mathbf{g})?\]\end{problem}

This minimal $\mathbf{f}$-weighted energy problem will be referred
to as the {\it Gauss variational problem\/}.
Cf.~\cite{DS,GR0,GR,NS,O,ST,Z4,Z5a,Z5,Z6}. Along with its
electrostatic interpretation, it has found various important
applications to approximation theory and to potential theory itself.

A minimizer~$\lambda$ is called an {\it equilibrium measure\/}
corresponding to the data~$\mathbf{A}$, $\mathbf{a}$, $\mathbf{g}$,
and~$\mathbf{f}$. The problem is said to be {\it solvable\/} if the
class $\mathcal G_{\mathbf{f}}(\mathbf{A},\mathbf{a},\mathbf{g})$ of
all those~$\lambda$ is nonempty.

\section{On uniqueness of equilibrium measures}

\begin{lemma}\label{lemma:unique}
If $\lambda$ and $\hat\lambda$ both belong to $\mathcal
G_{\mathbf{f}}(\mathbf{A},\mathbf{a},\mathbf{g})$, then\footnote{It
will also be shown below (see Corollary~\ref{cor:unique}) that
$\bigl\langle W_\lambda^i,\lambda^i\bigr\rangle= \bigl\langle
W_{\hat{\lambda}}^i,\hat{\lambda}^i\bigr\rangle$ for all $i\in I$.}

\begin{equation}\label{uniq1}
\|\lambda-\hat\lambda\|_{\mathcal E(\mathbf{A})}=0,\end{equation}
\begin{equation}\label{uniq3}
\langle\mathbf{f},\lambda\rangle=\langle\mathbf{f},\hat\lambda\rangle,\end{equation}
\begin{equation}\label{uniq2}
\mathbf{W}_\lambda(x)=\mathbf{W}_{\hat\lambda}(x)\quad\mbox{n.\,e.
in \ }\mathrm X.\end{equation}
\end{lemma}

\begin{proof}Since the class
$\mathcal E_{\mathbf{f}}(\mathbf{A},\mathbf{a},\mathbf{g})$ is
convex, we conclude from (\ref{G}), (\ref{wen}), and~(\ref{Ren})
that
\[4G_{\mathbf{f}}(\mathbf{A},\mathbf{a},\mathbf{g})\leqslant4G_{\mathbf{f}}\Bigl(\frac{\lambda+\hat{\lambda}}{2}\Bigr)=
\|R\lambda+R\hat{\lambda}\|^2+4\langle\mathbf{f},\lambda+\hat{\lambda}\rangle.\]
On the other hand, applying the parallelogram identity in the
pre-Hilbert space~$\mathcal E$ to $R\lambda$ and~$R\hat{\lambda}$
and then adding and subtracting
$4\langle\mathbf{f},\lambda+\hat{\lambda}\rangle$, we get
\[\|R\lambda-R\hat\lambda\|^2=
-\|R\lambda+R\hat{\lambda}\|^2-4\langle\mathbf{f},\lambda+\hat{\lambda}\rangle+2G_{\mathbf{f}}(\lambda)+
2G_{\mathbf{f}}(\hat{\lambda}).\] When combined with the preceding
relation, this yields
\[0\leqslant\|R\lambda-R\hat\lambda\|^2\leqslant-4G_{\mathbf{f}}(\mathbf{A},\mathbf{a},\mathbf{g})+2G_{\mathbf{f}}(\lambda)+
2G_{\mathbf{f}}(\hat{\lambda})=0,\] which establishes (\ref{uniq1})
because of~(\ref{seminorm}). In turn, (\ref{uniq1}) implies that
$\|\lambda\|^2=\|\hat{\lambda}\|^2$, whose subtraction from
$G_{\mathbf{f}}(\lambda)=G_{\mathbf{f}}(\hat{\lambda})$ results
in~(\ref{uniq3}). Due to Corollary~\ref{lemma:potequiv}, it can also
be concluded from~(\ref{uniq1}) that
$\kappa_\lambda(x)=\kappa_{\hat\lambda}(x)$ n.\,e. in~$\mathrm X$,
which together with~(\ref{wpot}) gives~(\ref{uniq2}).
\end{proof}

Thus, any two equilibrium measures (if exist) are equivalent in
$\mathcal E(\mathbf{A})$. Consequently, they are $R$-equivalent if
the kernel~$\kappa$ is strictly positive definite, and they are
equal if, moreover, all~$A_i$, $i\in I$, are mutually disjoint.

\section{Elementary properties of
$G_{\mathbf{f}}(\mathbf{A},\mathbf{a},\mathbf{g})$}

\subsection{Monotonicity and continuity of $G_{\mathbf{f}}(\,\cdot\,,\mathbf{a},\mathbf{g})$}
On $\mathfrak C=\mathfrak C(I^+,I^-)$, it is natural to introduce an
ordering relation~$\prec$ by declaring $\mathbf{A}'\prec\mathbf{A}$
to mean that $A_i'\subset A_i$ for all $i\in I$. Here,
$\mathbf{A}'=(A_i')_{i\in I}$. Then
$G_{\mathbf{f}}(\,\cdot\,,\mathbf{a},\mathbf{g})$ is a nonincreasing
function of a condenser, namely
\begin{equation}
G_{\mathbf{f}}(\mathbf{A},\mathbf{a},\mathbf{g})\leqslant
G_{\mathbf{f}}(\mathbf{A'},\mathbf{a},\mathbf{g})\quad\mbox{whenever
\ }\mathbf{A}'\prec\mathbf{A}. \label{increas'}
\end{equation}

Given $\mathbf{A}\in\mathfrak C$, we denote by
$\{\mathbf{K}\}_{\mathbf{A}}$ the increasing family of all compact
condensers $\mathbf{K}=(K_i)_{i\in I}\in\mathfrak C$ such that
$\mathbf{K}\prec\mathbf{A}$.

\begin{lemma}\label{lemma.cont} If $\mathbf{K}$ ranges over $\{\mathbf{K}\}_{\mathbf{A}}$, then
\begin{equation}
G_{\mathbf{f}}(\mathbf{A},\mathbf{a},\mathbf{g})=\lim_{\mathbf{K}\uparrow\mathbf{A}}\,
G_{\mathbf{f}}(\mathbf{K},\mathbf{a},\mathbf{g}).\label{cont}
\end{equation}\end{lemma}

\begin{proof} We can certainly assume that $G_{\mathbf{f}}(\mathbf{A},\mathbf{a},\mathbf{g})<\infty$,
since otherwise (\ref{cont})~follows at once from~(\ref{increas'}).
Then the set $\mathcal
E_{\mathbf{f}}(\mathbf{A},\mathbf{a},\mathbf{g})$ must be nonempty;
fix~$\mu$, one of its elements. Given
$\mathbf{K}\in\{\mathbf{K}\}_{\mathbf{A}}$ and $i\in I$, let
$\mu^i_{\mathbf{K}}$ denote the trace of~$\mu^i$ upon~$K_i$, i.\,e.,
$\mu^i_{\mathbf{K}}:=\mu_{K_i}^i$. Applying Lemma~1.2.2
from~\cite{F1} to~$g_i$, $f_{i1}$, $f_{i2}$, and~$\kappa$, we
conclude that
\begin{align}
\langle g_i,\mu^i\rangle&=\lim_{\mathbf{K}\uparrow\mathbf{A}}\,\langle g_i,\mu_{\mathbf{K}}^i\rangle,\qquad i\in I,\label{w}\\
\langle f_i,\mu^i\rangle&=\lim_{\mathbf{K}\uparrow\mathbf{A}}\,\langle f_i,\mu_{\mathbf{K}}^i\rangle,\qquad i\in I,\label{w'}\\
\kappa(\mu^i,\mu^j)&=\lim_{\mathbf{K}\uparrow\mathbf{A}}\,\kappa(\mu_{\mathbf{K}}^i,\mu_{\mathbf{K}}^j),\qquad
i,j\in I.\label{ww}
\end{align}

Fix $\varepsilon>0$. By~(\ref{w})--(\ref{ww}), for every $i\in I$
one can choose a compact set $K_i^0\subset A_i$ such that, for all
compact sets~$K_i$ with the properties $K_i^0\subset K_i\subset
A_i$, the following relations hold:
\begin{equation}
\frac{a_i}{\langle
g_i,\mu^i_{K_i}\rangle}<1+\varepsilon\,i^{-2},\label{unif2}
\end{equation}
\begin{equation}
\bigl|\langle f_i,\mu^i\rangle-\langle
f_i,\mu^i_{K_i}\rangle\bigr|<\varepsilon\,i^{-2},\label{unif2'}
\end{equation}
\begin{equation}
\bigl|\|\mu^i\|^2-\|\mu^i_{K_i}\|^2\bigr|<\varepsilon^2i^{-4}.\label{unif1}
\end{equation}
Having denoted $\mathbf{K}^0:=(K_i^0)_{i\in I}$, for every
$\mathbf{K}\in\{\mathbf{K}\}_{\mathbf{A}}$ that
follows~$\mathbf{K}^0$ we set
\begin{equation}\label{hatmu}\hat{\mu}^i_{\mathbf{K}}:=\frac{a_i}{\langle
g_i,\mu_{\mathbf{K}}^i\rangle}\,\mu_{\mathbf{K}}^i,\quad i\in
I.\end{equation} Then
$\hat{\mu}_{\mathbf{K}}:=\bigl(\hat{\mu}^i_{\mathbf{K}}\bigr)_{i\in
I}\in\mathcal E(\mathbf{K},\mathbf{a},\mathbf{g})$, the finiteness
of the energy being obtained from~(\ref{unif1}) and
Lemma~\ref{enfinite}. Furthermore, since $\sum_{i\in I}\,\langle
f_i,\mu^i\rangle$ is absolutely convergent, so is $\sum_{i\in
I}\,\langle f_i,\hat{\mu}^i_{\mathbf{K}}\rangle$, which is clear
from~(\ref{unif2}) and~(\ref{unif2'}). Therefore actually
$\hat{\mu}_{\mathbf{K}}\in\mathcal
E_{\mathbf{f}}(\mathbf{K},\mathbf{a},\mathbf{g})$, and consequently
\begin{equation}
G_{\mathbf{f}}(\hat{\mu}_{\mathbf{K}})\geqslant
G_{\mathbf{f}}(\mathbf{K},\mathbf{a},\mathbf{g}).\label{www}\end{equation}

We next proceed by showing that
\begin{equation}
G_{\mathbf{f}}(\mu)=\lim_{\mathbf{K}\uparrow\mathbf{A}}\,G_{\mathbf{f}}(\hat{\mu}_{\mathbf{K}}).\label{4w}\end{equation}
To this end, it can be assumed that $\kappa\geqslant0$; for if not,
then $\mathbf{A}$ must be finite since $\mathrm X$ is compact, and
(\ref{4w}) follows from~(\ref{w})--(\ref{ww}). Therefore, for all
$\mathbf{K}\succ\mathbf{K}_0$ and $i\in I$ we get
\begin{equation}
\|\mu^i_{\mathbf{K}}\|\leqslant\|\mu^i\|\leqslant\|R\mu^++R\mu^-\|,
\label{unif3}
\end{equation}
\begin{equation}
\|\mu^i-\mu^i_{\mathbf{K}}\|<\varepsilon\,i^{-2},\label{uniff}
\end{equation}
the latter being clear from (\ref{unif1}) because of
$\kappa(\mu^i_{\mathbf{K}},\mu^i-\mu^i_{\mathbf{K}})\geqslant0$.
Also observe that
\begin{equation*}
\begin{split}
\bigl|&\|\mu\|^2-\|\hat{\mu}_{\mathbf{K}}\|^2\bigr|\leqslant\sum_{i,j\in
I}\,\Bigl|\kappa(\mu^i,\mu^j)- \frac{a_i}{\langle
g_i,\mu_{\mathbf{K}}^i\rangle}\frac{a_j}{\langle
g_j,\mu_{\mathbf{K}}^j\rangle}\,
\kappa(\mu_{\mathbf{K}}^i,\mu_{\mathbf{K}}^j)\Bigr|\\[3pt]
&{}\leqslant\sum_{i,j\in
I}\,\Bigl[\kappa(\mu^i-\mu^i_{\mathbf{K}},\mu^j)+\kappa(\mu^i_{\mathbf{K}},\mu^j-\mu^j_{\mathbf{K}})+
\Bigl(\frac{a_i}{\langle
g_i,\mu_{\mathbf{K}}^i\rangle}\frac{a_j}{\langle
g_j,\mu_{\mathbf{K}}^j\rangle}-1\Bigr)\,\kappa(\mu^i_{\mathbf{K}},\mu^j_{\mathbf{K}})\Bigr].
\end{split}
\end{equation*}
When combined with (\ref{unif2}), (\ref{unif2'}), (\ref{unif3}),
and~(\ref{uniff}), this yields
\[
\bigl|G_{\mathbf{f}}(\mu)-G_{\mathbf{f}}(\hat{\mu}_{\mathbf{K}})\bigr|\leqslant
M\varepsilon\quad\mbox{for all \ }\mathbf{K}\succ\mathbf{K}_0,\]
where $M$ is finite and independent of~$\mathbf{K}$, and the
required relation~(\ref{4w}) follows.

Substituting (\ref{www}) into~(\ref{4w}), in view of the arbitrary
choice of $\mu\in\mathcal
E_{\mathbf{f}}(\mathbf{A},\mathbf{a},\mathbf{g})$ we get
\[
G_{\mathbf{f}}(\mathbf{A},\mathbf{a},\mathbf{g})\geqslant
\lim_{\mathbf{K}\uparrow\mathbf{A}}\,G_{\mathbf{f}}(\mathbf{K},\mathbf{a},\mathbf{g})\|^2.
\]
Since the converse inequality is obvious from~(\ref{increas'}), the
proof is complete.
\end{proof}

Let $\mathcal E_{\mathbf{f}}^0(\mathbf{A},\mathbf{a},\mathbf{g})$
denote the class of all $\mu\in\mathcal
E_{\mathbf{f}}(\mathbf{A},\mathbf{a},\mathbf{g})$ such that, for
every $i\in I$, the support~$S(\mu^i)$ of~$\mu^i$ is compact.

\begin{corollary} \label{cor:compact} The value $G_{\mathbf{f}}(\mathbf{A},\mathbf{a},\mathbf{g})$
remains unchanged if the class $\mathcal
E_{\mathbf{f}}(\mathbf{A},\mathbf{a},\mathbf{g})$ in its definition
is replaced by $\mathcal
E_{\mathbf{f}}^0(\mathbf{A},\mathbf{a},\mathbf{g})$. That is,
\[G_{\mathbf{f}}(\mathbf{A},\mathbf{a},\mathbf{g})=\inf_{\mu\in\mathcal E_{\mathbf{f}}^0(\mathbf{A},\mathbf{a},\mathbf{g})}\,
G_{\mathbf{f}}(\mu).\]\end{corollary}

\subsection{When does $G_{\mathbf{f}}(\mathbf{A},\mathbf{a},\mathbf{g})<\infty$
hold?}\label{nonz} Let $C(E)$ denote the interior capacity of a
set~$E\subset\mathrm X$. Given $\mathbf{g}=(g_i)_{i\in I}$, we also
write
\[g_{i,\inf}:=\inf_{x\in A_i}\,g_i(x),\qquad g_{i,\sup}:=\sup_{x\in A_i}\,g_i(x).\]

This section provides necessary and (or) sufficient conditions for
the class $\mathcal
E_{\mathbf{f}}(\mathbf{A},\mathbf{a},\mathbf{g})$ to be nonempty or,
which is equivalent, for
\begin{equation}
G_{\mathbf{f}}(\mathbf{A},\mathbf{a},\mathbf{g})<\infty.\label{nonzero1}
\end{equation}

\begin{lemma}\label{lemma:nonzero}For {\rm(\ref{nonzero1})} to hold, it is
necessary that
\begin{equation}\label{nec}C\bigl(\{x\in A_i: \ |f_i(x)|<\infty\}\bigr)\ne0\quad\mbox{for all \ }i\in I.\end{equation}
If $\mathbf{A}$ is finite, then {\rm(\ref{nonzero1})} and
{\rm(\ref{nec})} are actually equivalent.\end{lemma}

\begin{proof}If {\rm(\ref{nonzero1})} holds, then by Corollary~\ref{cor:compact}
there is $\mu\in\mathcal
E^0_{\mathbf{f}}(\mathbf{A},\mathbf{a},\mathbf{g})$. Assume, on the
contrary, that $C\bigl(\{x\in A_{i_0}:\
|f_{i_0}(x)|<\infty\}\bigr)=0$ for some $i_0\in I$. Since
$\mu^{i_0}$ has finite energy and is compactly supported
in~$A_{i_0}$, \cite[Lemma~2.3.1]{F1} yields that
$|f_{i_0}(x)|=\infty$ $\mu^{i_0}$-al\-most everywhere
($\mu^{i_0}$-a.\,e.) in~$\mathrm X$. This is impossible, for
$\mu^{i_0}$ is nonzero while $\langle\mathbf{f},\mu\rangle$ is
finite.

Assuming now $\mathbf{A}$ to be finite, we proceed by proving that
(\ref{nec}) implies (\ref{nonzero1}). For each $i\in I$, the set
$E_i:=\{x\in A_i:\ |f_i(x)|<\infty\}$ can be written as the union of
$E^n_i$, $n\in\mathbb N$, where $E^n_i:=\{x\in A_i: \
|f_i(x)|\leqslant n\}$. Taking into account that $E^n_i$ are
increasing and universally meas\-ur\-able, from \cite[Lemma
2.3.3]{F1} we get $C(E_i)=\lim_{n\to\infty}\,C(E_i^n)$. Since
$C(E_i)>0$ while $\mathbf{A}$ is finite, one can choose~$n_0$ so
that $C(E_i^{n_0})>0$ for all $i\in I$. Consequently, for every
$i\in I$ there is a probability measure~$\omega_i$ of finite energy,
compactly supported in~$E_i^{n_0}$.

The function $g_i$, being continuous, is bounded on $S(\omega_i)$;
hence $0<\langle g_i,\omega_i\rangle<\infty$. Writing
\[\hat{\omega}^i:=\frac{a_i\omega_i}{\langle g_i,\omega_i\rangle},\quad i\in I,\]
we obtain $\hat{\omega}:=(\hat{\omega}^i)_{i\in I}\in\mathcal
E(\mathbf{A},\mathbf{a},\mathbf{g})$. Since $|\langle
f_i,\hat{\omega}^i\rangle|\leqslant n_0\hat{\omega}^i(\mathrm
X)<\infty$ for all $i\in I$, we actually have
$\hat{\omega}\in\mathcal
E_{\mathbf{f}}(\mathbf{A},\mathbf{a},\mathbf{g})$, and
the desired relation (\ref{nonzero1}) follows.
\end{proof}

However, if $\mathbf{A}$ is infinite, assuming only (\ref{nec}) is
not enough to guarantee (\ref{nonzero1}); then sufficient conditions
for (\ref{nonzero1}) to hold can be formulated as follows.

\begin{lemma} Assume there are constants~$M<\infty$ and $\delta>0$, both independent
of~$i$, such that
\begin{equation*}
C\bigl(\{x\in A_i: \ |f_i(x)|\leqslant
M\}\bigr)>\delta\quad\mbox{for all \ }i\in I.
\end{equation*}
Then {\rm(\ref{nonzero1})} is true whenever
\begin{equation}\label{s}\sum_{i\in
I}\,a_ig_{i,\inf}^{-1}<\infty.\end{equation}
\end{lemma}

\begin{proof}For every $i\in I$, we denote $E_i^M:=\{x\in A_i: \ |f_i(x)|\leqslant
M\}$ and choose a probability measure $\omega_i\in\mathcal
E^+(E_i^M)$ so that \[\|\omega_i\|^2\leqslant
C(E_i^M)^{-1}+\delta<\delta+\delta^{-1}.\] Defining
$\hat{\omega}^i$, $i\in I$, by the same formula as in the preceding
proof, we then obtain, by~(\ref{s}), \[\sum_{i\in
I}\,\|\hat{\omega}^i\|^2\leqslant\bigl[\delta+\delta^{-1}\bigr]\sum_{i\in
I}\,a^2_ig_{i,\inf}^{-2}<\infty\] and hence, by
Lemma~\ref{enfinite}, $\hat{\omega}:=(\hat{\omega}^i)_{i\in
I}\in\mathcal E(\mathbf{A},\mathbf{a},\mathbf{g})$. Since,
by~(\ref{s}),
\[\sum_{i\in
I}\,|\langle f_i,\hat{\omega}^i\rangle|\leqslant M\sum_{i\in
I}\,\hat{\omega}^i(\mathrm X)\leqslant M\sum_{i\in
I}\,a_ig_{i,\inf}^{-1}<\infty,\] we actually have
$\hat{\omega}\in\mathcal
E_{\mathbf{f}}(\mathbf{A},\mathbf{a},\mathbf{g})$, and
the claimed conclusion follows.
\end{proof}

\section{Description of the $\mathbf{f}$-weighted equilibrium
potentials}\label{sec:desc}

Given a set $E\subset\mathrm X$ of interior capacity nonzero and a
universally measurable function~$\psi$ bounded from below nearly
everywhere in~$E$, write
$$
"\!\inf_{x\in E}\!"\,\,\psi(x):=\sup\,\bigl\{q: \ \psi(x)\geqslant
q\quad\mbox{n.\,e.~in \ } E\bigr\}.
$$
Then \[ \psi(x)\geqslant"\!\inf_{x\in
E}\!"\,\,\psi(x)\quad\mbox{n.\,e.~in \ } E,\] which follows from the
fact that the union of a sequence of sets $U_n\cap E$ with
$C(U_n\cap E)=0$ is of interior capacity zero as well, provided
$U_n$, $n\in\mathbb N$, are universally measurable whereas $E$ is
arbitrary (see the corollary to~Lemma~2.3.5 in~\cite{F1} and the
remark attached to it).

\subsection{Variational inequalities for the $\mathbf{f}$-weighted equilibrium potentials}

Throughout Sect.~\ref{sec:desc} we assume that an equilibrium
measure $\lambda$ exists (see~Theorem~\ref{exist} for conditions
ensuring the solvability of the Gauss variational problem). Then,
for every $i\in I$, $W_\lambda^i(x)$ is defined and ${}\ne-\infty$
n.\,e.~in~$A_i$, while $C(A_i)>0$ as a consequence
of~Lemma~\ref{lemma:nonzero}.

\begin{theorem}\label{th:descpot1}
For all $\lambda\in\mathcal
G_{\mathbf{f}}(\mathbf{A},\mathbf{a},\mathbf{g})$ and $i\in I$,
\begin{equation}\label{descpot1}
a_i\,W_\lambda^i(x)\geqslant\bigl\langle
W_\lambda^i,\lambda^i\bigr\rangle\,g_i(x)\quad\mbox{n.\,e. in \ }
A_i.
\end{equation}
\end{theorem}

\begin{proof}Indeed, $\lambda^i$ is a solution to the problem of minimizing
$G_{\tilde{f_i}}(\nu)=\|\nu\|^2+2\langle\tilde{f_i},\nu\rangle$,
where
\begin{equation*}\label{th:desc1:1}\tilde{f_i}(x):=f_i(x)+\alpha_i\sum_{j\in
I,\,j\ne i}\,\alpha_j\kappa(x,\lambda^j)\end{equation*} and $\nu$
ranges over the class $\mathcal E_{\tilde{f_i}}(A_i,a_i,g_i)$.
Applying \cite[Th.~2.1]{O}, we arrive
at~(\ref{descpot1}).
\end{proof}

In the following assertion we additionally assume that, for each
$i\in I$, either $g_{i,\inf}>0$ or $A_i$~can be written as a
countable union of compact sets. Then every~$A_i$ is a countable
union of $\nu^i$-integrable sets, where $\nu\in\mathfrak
M(\mathbf{A},\mathbf{a},\mathbf{g})$ is arbitrarily given, and hence
any locally $\nu^i$-negligible subset of~$A_i$ is
$\nu^i$-negligible.

\begin{corollary}\label{cor:descpot1} For all $\lambda\in\mathcal G_{\mathbf{f}}(\mathbf{A},\mathbf{a},\mathbf{g})$
and $i\in I$,
\begin{equation}\label{descpot1''}
a_i\,W_\lambda^i(x)=\bigl\langle
W_\lambda^i,\lambda^i\bigr\rangle\,g_i(x)\quad\lambda^i\mbox{-a.\,e.
in \ } \mathrm X.
\end{equation}
\end{corollary}

\begin{proof} Since $\lambda^i$ has finite energy, the set of all $x\in A_i$ for which the inequality
in~(\ref{descpot1}) fails to hold is locally $\lambda^i$-negligible
by~\cite[Lemma~2.3.1]{F1} and, hence, it is $\lambda^i$-negligible
(cf.~the note followed by the corollary). Hence, (\ref{descpot1''})
must be true, for if not, we would arrive at a contradiction by
integrating the inequality in~(\ref{descpot1}) with respect
to~$\lambda^i$.
\end{proof}

\begin{theorem}\label{th:descpot2}
Assume $\kappa$ is continuous on $A^+\times A^-$ and satisfies the
condition
\begin{equation}\label{Kbounded}
\sup_{x\in K,\,y\in A^-}\,\kappa(x,y)<\infty\quad\mbox{for all
compact \ }K\subset A^+
\end{equation}
and that obtained from {\rm(\ref{Kbounded})} when the indices $+$
and~$-$ are reversed. Let moreover $f_i\in\mathrm\Phi(\mathrm X)$
for all $i\in I$, and let {\rm(\ref{s})} hold true. For every
$\lambda\in\mathcal
G_{\mathbf{f}}(\mathbf{A},\mathbf{a},\mathbf{g})$, then
\begin{equation}\label{descpot2'}
a_i\,W_\lambda^i(x)\leqslant\bigl\langle
W_\lambda^i,\lambda^i\bigr\rangle\,g_i(x)\quad\mbox{for all \ } x\in
S(\lambda^i)
\end{equation}
and, hence,
\begin{equation}\label{descpot2''}
a_i\,W_\lambda^i(x)=\bigl\langle
W_\lambda^i,\lambda^i\bigr\rangle\,g_i(x)\quad\mbox{n.\,e. in \ }
S(\lambda^i).
\end{equation}
\end{theorem}

\begin{proof}Fix $i\in I$ (say $i\in I^+$). We begin by verifying that $W_\mu^i$, where
$\mu\in\mathcal E_{\mathbf{f}}(\mathbf{A},\mathbf{a},\mathbf{g})$ is
given, is lower semicontinuous on~$A_i$. To this end, it is enough
to show that so is~$-\kappa(\,\cdot\,,R\mu^-)$.

Having fixed a point $x_0\in A_i$ and its compact neighborhood
$V_{x_0}\subset A_i$, let us consider a function $\kappa^*(x,y)$ on
$V_{x_0}\times A^-$, defined by the formula
\begin{equation}\label{kappastar}
\kappa^*(x,y):=-\kappa(x,y)+\sup_{{x'}\in V_{x_0},\,{y'}\in
A^-}\,\kappa({x'},{y'}).
\end{equation}
Under the assumptions of the theorem,  $\kappa^*$ is nonnegative and
continuous; hence,
\[\kappa^*(x,R\mu^-)=\int\kappa^*(x,y)\,dR\mu^-(y),\quad x\in
V_{x_0},\] being the potential of the nonnegative measure~$R\mu^-$
with respect to the kernel~$\kappa^*$, is lower semicontinuous.

On the other hand, it follows from~(\ref{s}) that $R\mu^-$ is
bounded. Integrating~(\ref{kappastar}) with respect to~$R\mu^-$, we
conclude from~(\ref{Kbounded}) that $\kappa^*(x,R\mu^-)$, $x\in
V_{x_0}$, coincides up to a finite summand with the restriction
of~$-\kappa(x,R\mu^-)$ to~$V_{x_0}$. What has been shown just above
therefore implies that $-\kappa(\,\cdot\,,R\mu^-)$ is lower
semicontinuous on~$A_i$. Hence, so is $W_\mu^i$.

To complete the proof, fix $\lambda\in\mathcal
G_{\mathbf{f}}(\mathbf{A},\mathbf{a},\mathbf{g})$ and $x\in
S(\lambda^i)$, and let $\mathcal B(x)$ be the family of all
neighborhoods of~$x$ in~$A_i$, directed by~${}\subset{}$. For every
$U\in\mathcal B(x)$, we have $\lambda^i(U)>0$; hence,
by~(\ref{descpot1''}), one can choose a point~$x_U\in U$ so that
\[a_i\,W_\lambda^i(x_U)=\bigl\langle
W_\lambda^i,\lambda^i\bigr\rangle\,g_i(x_U).\] Since the net
$\bigl(x_U\bigr)_{U\in\mathcal B(x)}$ converges to~$x$, this
proves~(\ref{descpot2'}) because $W_\lambda^i$ is lower
semicontinuous on~$A_i$ while $g_i$ is continuous. Finally,
combining (\ref{descpot1}) and~(\ref{descpot2'})
gives~(\ref{descpot2''}).
\end{proof}

\subsection{Characteristic properties of equilibrium
measures}\label{sec:char} Observing that
\begin{equation}\label{star}G_{\mathbf{f}}(\mathbf{A},\mathbf{a},\mathbf{g})=G_{\mathbf{f}}(\lambda)=
\sum_{i\in I}\,\bigl\langle
W_\lambda^i,\lambda^i\bigr\rangle+\langle\mathbf{f},\lambda\rangle,\end{equation}
we proceed by showing that (\ref{descpot1}), (\ref{descpot1''})
and~(\ref{star}) serve as characteristic properties of~$\lambda$.

\begin{theorem}\label{th:char}
Given $\mu\in\mathcal
E_{\mathbf{f}}(\mathbf{A},\mathbf{a},\mathbf{g})$, suppose there are
numbers $\eta_i$ such that, for all $i\in I$, either
{\rm(\ref{descpot1'})} and~{\rm(\ref{descpot2})} or
{\rm(\ref{descpot1,'})} and~{\rm(\ref{descpot2,,})} hold true, where
\begin{align}\label{descpot1'}
a_i\,W_\mu^i(x)&\geqslant\eta_i\,g(x)\quad\mbox{n.\,e. in \ } A_i,\\
\label{descpot2} G_{\mathbf{f}}(\mu)&\leqslant\sum_{i\in
I}\,\eta_i+\langle\mathbf{f},\mu\rangle
\end{align}
and
\begin{align}\label{descpot1,'}
a_i\,W_\mu^i(x)&\leqslant\eta_i\,g(x)\quad\mu^i\mbox{-a.\,e. in \ } \mathrm X,\\
\label{descpot2,,}
G_{\mathbf{f}}(\mathbf{A},\mathbf{a},\mathbf{g})&\geqslant\sum_{i\in
I}\,\eta_i+\langle\mathbf{f},\mu\rangle.
\end{align}
Then $\mu$ belongs to $\mathcal
G_{\mathbf{f}}(\mathbf{A},\mathbf{a},\mathbf{g})$ and
\begin{equation}\label{eta}\eta_i=\bigl\langle
W_\mu^i,\mu^i\bigr\rangle\quad\mbox{for all \ } i\in
I.\end{equation}
\end{theorem}

\begin{proof}Assuming (\ref{descpot1'}) and~(\ref{descpot2}) to
hold, fix $\nu\in\mathcal
E^0_{\mathbf{f}}(\mathbf{A},\mathbf{a},\mathbf{g})$. Since $\nu^i$
is of finite energy and compactly supported in~$A_i$,
\cite[Lemma~2.3.1]{F1} shows that the inequality
in~(\ref{descpot1'}) holds $\nu^i$-a.\,e.~in~$\mathrm X$. This gives
\begin{equation}\label{prdesc}\bigl\langle
W_\mu^i,\nu^i\bigr\rangle\geqslant\eta_i\quad\mbox{for all \ } i\in
I.\end{equation} Summing up these inequalities and then
substituting~(\ref{descpot2}) into the result obtained, we get
\[\kappa(\nu,\mu)+\langle\mathbf{f},\nu\rangle\geqslant\|\mu\|^2+\langle\mathbf{f},\mu\rangle,\]
which in turn yields
\[G_{\mathbf{f}}(\nu)-G_{\mathbf{f}}(\mu)\geqslant\|\nu-\mu\|^2.\]
Application of Corollary~\ref{cor:compact} therefore implies that
$\mu$ is an equilibrium measure.

Further, for all $\mathbf{K}\in\{\mathbf{K}\}_{\mathbf{A}}$ large
enough consider $\hat{\mu}^i_{\mathbf{K}}$ defined by~(\ref{hatmu}).
Applying~(\ref{prdesc}) to~$\hat{\mu}^i_{\mathbf{K}}$ instead
of~$\nu^i$ and then letting $\mathbf{K}\uparrow\mathbf{A}$, by
arguments similar to those used in the proof of
Lemma~\ref{lemma.cont} we get $\bigl\langle
W^i_\mu,\mu^i\bigr\rangle\geqslant\eta_i$ for all $i\in I$. Summing
up these inequalities  and then comparing the result obtained with
(\ref{star}) for~$\lambda$ replaced by~$\mu$ and (\ref{descpot2}),
we obtain~(\ref{eta}).

Since the remaining case can be handled in a similar way, the proof
is complete.
\end{proof}

\begin{corollary}\label{cor:unique}
$\bigl\langle W_\lambda^i,\lambda^i\bigr\rangle=\bigl\langle
W_{\hat{\lambda}}^i,\hat{\lambda}^i\bigr\rangle$ for any
$\lambda,\,\hat{\lambda}\in\mathcal
G_{\mathbf{f}}(\mathbf{A},\mathbf{a},\mathbf{g})$ and all $i\in
I$.\end{corollary}
\begin{corollary}Given $\lambda\in\mathcal
G_{\mathbf{f}}(\mathbf{A},\mathbf{a},\mathbf{g})$, we have
\begin{equation}\bigl\langle
W_\lambda^i,\lambda^i\bigr\rangle="\!\inf_{x\in
A_i}\!"\,\,\frac{a_i\,W_\lambda^i(x)}{g(x)}\quad\mbox{for all \ }
i\in I\label{essinf}\end{equation} and, hence,
\[G_{\mathbf{f}}(\lambda)=\sum_{i\in I}\,"\!\inf_{x\in
A_i}\!"\,\,\frac{a_i\,W_\lambda^i(x)}{g(x)}+\langle\mathbf{f},\lambda\rangle.\]
\end{corollary}

\subsection{$\mathbf{f}$-weighted equilibrium
constants}

\begin{definition}\label{eqcont} We shall call the numbers
$\bigl\langle W_\lambda^i,\lambda^i\bigr\rangle$, $i\in I$, where
$\lambda\in\mathcal
G_{\mathbf{f}}(\mathbf{A},\mathbf{a},\mathbf{g})$ is arbitrarily
given, the $\mathbf{f}$-{\it weighted equilibrium constants\/}
corresponding to the data $\mathbf{A}$, $\mathbf{a}$, $\mathbf{g}$,
and~$\mathbf{f}$.\end{definition}

These constants do not depend on the choice of $\lambda\in\mathcal
G_{\mathbf{f}}(\mathbf{A},\mathbf{a},\mathbf{g})$, which is clear
from Corollary~\ref{cor:unique}. They can also be uniquely
determined as $\eta_i$, $i\in I$, satisfying both the
relations~(\ref{descpot1'}) and~(\ref{descpot2})
with~$\lambda\in\mathcal
G_{\mathbf{f}}(\mathbf{A},\mathbf{a},\mathbf{g})$ in place of~$\mu$.
Another alternative definition of the $\mathbf{f}$-weighted
equilibrium constants can be given by~(\ref{essinf}).

\section{Equilibrium measures: existence and $\mathbf{A}$-vague compactness.
 Statements on continuity}

Assume for a moment that a condenser $\mathbf{A}$~is compact. Then
the class $\mathfrak M(\mathbf{A},\mathbf{a},\mathbf{g})$ is
$\mathbf{A}$-vaguely bounded and closed and hence, by
Lemma~\ref{lem:vaguecomp}, it is $\mathbf{A}$-vaguely compact. If
moreover $\mathbf{A}$~is finite, $\kappa$ is continuous
on~$A^+\times A^-$, while $f_i\in\mathrm\Phi(\mathrm X)$ for all
$i\in I$, then $G_{\mathbf{f}}(\mu)$ is $\mathbf{A}$-va\-gue\-ly
lower semicontinuous on~$\mathcal E(\mathbf{A})$ and, therefore, the
existence of equilibrium measures~$\lambda$ immediately follows.
See~\cite[Th.~2.30]{O}; cf.~also~\cite{GR0,GR,NS,ST}.

However, these arguments break down if any of the above assumptions
is dropped. In particular, $\mathfrak
M(\mathbf{A},\mathbf{a},\mathbf{g})$ is no longer
$\mathbf{A}$-vaguely compact if $\mathbf{A}$ is noncompact.

To solve the problem on the existence of equilibrium measures in the
general case where a condenser $\mathbf{A}$ is infinite and (or)
noncompact, we develop an approach based on both the
$\mathbf{A}$-vague and strong topologies in the semimetric
space~$\mathcal E(\mathbf{A})$, introduced for measures of finite
dimensions in~\cite{Z4,Z5a,Z5,Z6}.

\subsection{Standing assumptions}\label{sec:standing}

Unless explicitly stated otherwise, in all that follows it is
required that the kernel~$\kappa$ is consistent and either
$I^-=\varnothing$, or (\ref{s}) and the following condition are both
satisfied:
\begin{equation}\sup_{x\in A^+,\ y\in
A^-}\,\kappa(x,y)<\infty.\label{bou}
\end{equation}
It will also be assumed that
$G_{\mathbf{f}}(\mathbf{A},\mathbf{a},\mathbf{g})<\infty$, which
certainly involves no loss of generality, since otherwise the Gauss
variational problem makes no sense; see Sect.~\ref{nonz} for
necessary and (or) sufficient conditions for this to hold.

Throughout Sections~\ref{sec:ex} and~\ref{appr} we shall also
suppose one of the following Cases~I, II, or~III to occur:
\begin{itemize}
\item[\rm I.] There exists a vector measure $\nu\in\mathcal E(\mathbf{A})$ such that
$\mathbf{f}=\kappa_\nu$;
\item[\rm II.] There exists $\sigma\in\mathcal E$ such that
$f_i=\alpha_i\kappa(\,\cdot\,,\sigma)$ for all $i\in I$;
\item[\rm III.] $f_i\in\mathrm\Phi(\mathrm X)$ for all $i\in I$.
\end{itemize}

\begin{remark} In all the Cases I, II, or III, the restrictions on~$\mathbf{f}$ that have
been imposed in~Sect.~\ref{sec:statement} do hold
automatically.\end{remark}

\begin{remark}
Note that the above assumptions on a kernel are not too restrictive.
In particular, they all are satisfied by the Newtonian, Riesz, or
Green kernels in~$\mathbb R^n$, $n\geqslant2$, provided the
Euclidean distance between $A^+$ and $A^-$ is nonzero.\end{remark}

\subsection{Statements on existence and $\mathbf{A}$-vague compactness}\label{sec:ex}
\begin{theorem}\label{exist}
Under the standing assumptions, let moreover for every $i\in I$
either $g_{i,\sup}<\infty$ or there exist $r_i\in(1,\infty)$ and
$\omega_i\in\mathcal E$ such that
\begin{equation}
g_i^{r_i}(x)\leqslant\kappa(x,\omega_i)\quad\mbox{n.\,e. in \ } A_i.
\label{growth}
\end{equation}
If, in addition, $A_i$ either is compact or has finite interior
capacity\footnote{Note that a compact set $K\subset\mathrm X$ might
be of infinite capacity; $C(K)$ is necessarily finite provided the
kernel is strictly positive definite~\cite{F1}. On the other hand,
even for the Newtonian kernel sets of finite capacity might be
noncompact (see~\cite{L}).}, then the class of equilibrium measures
$\mathcal G_{\mathbf{f}}(\mathbf{A},\mathbf{a},\mathbf{g})$ is
nonempty and $\mathbf{A}$-vaguely compact.
\end{theorem}

\begin{corollary}\label{cor:exist} If $\mathbf{A}=\mathbf{K}$ is compact, then $\mathcal
G_{\mathbf{f}}(\mathbf{A},\mathbf{a},\mathbf{g})$ is nonempty and
$\mathbf{A}$-vaguely compact.
\end{corollary}
\begin{proof}This is an immediate consequence of Theorem~\ref{exist},
for $g_i$ is bounded on~$K_i$.\end{proof}

\subsection{On continuity of equilibrium measures and $\mathbf{f}$-weighted equilibrium constants}\label{appr}
When approaching $\mathbf{A}$ by the increasing
family~$\{\mathbf{K}\}_{\mathbf{A}}$ of the compact condensers
$\mathbf{K}\prec\mathbf{A}$, we shall always suppose all
those~$\mathbf{K}$ to satisfy the assumption
$G_{\mathbf{f}}(\mathbf{K},\mathbf{a},\mathbf{g})<\infty$. This
involves no loss of generality, which is clear from the
assumption~(\ref{nonzero1}) and Lemma~\ref{lemma.cont}. Choose an
equilibrium measure $\lambda_{\mathbf{K}}\in\mathcal
G_{\mathbf{f}}(\mathbf{K},\mathbf{a},\mathbf{g})$~--- its existence
has been ensured by Corollary~\ref{cor:exist}.

\begin{theorem}\label{cor:cont} Let all
the conditions of Theorem~{\rm\ref{exist}} be satisfied. Then every
$\mathbf{A}$-vague cluster point of
$(\lambda_{\mathbf{K}})_{\mathbf{K}\in\{\mathbf{K}\}_{\mathbf{A}}}$
(such a cluster point exists) belongs to $\mathcal
G_{\mathbf{f}}(\mathbf{A},\mathbf{a},\mathbf{g})$. Furthermore, if
$\lambda_{\mathbf{A}}\in\mathcal
G_{\mathbf{f}}(\mathbf{A},\mathbf{a},\mathbf{g})$ is arbitrarily
given, then
\begin{equation}\label{cor:cont1}\lim_{\mathbf{K}\uparrow\mathbf{A}}\,\|\lambda_{\mathbf{K}}-\lambda_{\mathbf{A}}\|^2=0,
\end{equation}
\begin{equation}\label{cor:cont2}\lim_{\mathbf{K}\uparrow\mathbf{A}}\,\langle\mathbf{f},\lambda_{\mathbf{K}}\rangle=
\langle\mathbf{f},\lambda_{\mathbf{A}}\rangle,\end{equation}
\begin{equation}\label{cor:cont3}\lim_{\mathbf{K}\uparrow\mathbf{A}}\,
\bigl\langle
W^i_{\lambda_{\mathbf{K}}},\lambda^i_{\mathbf{K}}\bigr\rangle=\bigl\langle
W^i_{\lambda_{\mathbf{A}}},\lambda^i_{\mathbf{A}}\bigr\rangle\quad\mbox{for
all \ }i\in I.\end{equation}
\end{theorem}

Thus, under the assumptions of Theorem~\ref{cor:cont}, if moreover
$\kappa$ is strictly positive definite and all~$A_i$, $i\in I$, are
mutually disjoint, then the (unique) equilibrium
measure~$\lambda_{\mathbf{K}}$ on~$\mathbf{K}$ converges both
$\mathbf{A}$-vaguely and strongly to the (unique) equilibrium
measure~$\lambda_{\mathbf{A}}$ on~$\mathbf{A}$.

The proofs of Theorems~\ref{exist} and~\ref{cor:cont}, to be given
in Sections~\ref{sec:proof.th.str} and~\ref{sec:proof.th.cont} below
(see also Sect.~\ref{sec:extremal} for auxiliary notions and
results), are based on a theorem on the strong completeness of
proper subspaces of the semimetric space~$\mathcal E(\mathbf{A})$,
which is a subject of the next section.

\section{Strong completeness of vector measures}\label{sec:11}

As always, assume all the standing assumptions, stated
in~Sect.~\ref{sec:standing}, to hold. Having denoted
\begin{equation*}
\mathfrak
M(\mathbf{A},\leqslant\!\mathbf{a},\mathbf{g}):=\bigl\{\mu\in\mathfrak
M(\mathbf{A}):\quad\langle g_i,\mu^i\rangle\leqslant a_i\mbox{ \ for
all \ } i\in I\bigr\},
\end{equation*}
we consider $\mathcal
E(\mathbf{A},\leqslant\!\mathbf{a},\mathbf{g}):=\mathfrak
M(\mathbf{A},\leqslant\!\mathbf{a},\mathbf{g})\cap\mathcal
E(\mathbf{A})$ to be a topological subspace of the semimetric space
$\mathcal E(\mathbf{A})$; the induced topology is likewise called
the {\it strong\/} topology.

Our purpose is to show that $\mathcal
E(\mathbf{A},\leqslant\!\mathbf{a},\mathbf{g})$ is strongly
complete.

\subsection{Auxiliary assertions}

\begin{lemma}\label{Mbounded} The class $\mathfrak M(\mathbf{A},\leqslant\!\mathbf{a},\mathbf{g})$ is
$\mathbf{A}$-vaguely bounded and, hence, $\mathbf{A}$-vaguely
compact.\end{lemma}

\begin{proof} Fix $i\in I$, and let a compact set $K_i\subset A_i$ be given. Since
$g_i$ is positive and continuous, the relation
\[
a_i\geqslant\langle g_i,\mu^i\rangle\geqslant\mu^i(K_i)\,\min_{x\in
K_i}\ g_i(x),\quad\mbox{where \ }\mu\in\mathfrak
M(\mathbf{A},\leqslant\!\mathbf{a},\mathbf{g}),\] yields
\[\sup_{\mu\in\mathfrak M(\mathbf{A},\leqslant\mathbf{a},\mathbf{g})}\,\mu^i(K_i)<\infty.\]
This implies that $\mathfrak
M(\mathbf{A},\leqslant\!\mathbf{a},\mathbf{g})$ is
$\mathbf{A}$-vaguely bounded and hence, by
Lemma~\ref{lem:vaguecomp}, $\mathbf{A}$-va\-guely relatively
compact. Since it is obviously $\mathbf{A}$-vaguely closed, the
lemma
follows.
\end{proof}

\begin{lemma}\label{lem:aux2} If a net
$(\mu_s)_{s\in S}\subset\mathcal
E(\mathbf{A},\leqslant\!\mathbf{a},\mathbf{g})$ is strongly bounded,
then its $\mathbf{A}$-vague cluster set is contained in $\mathcal
E(\mathbf{A},\leqslant\!\mathbf{a},\mathbf{g})$.\end{lemma}

\begin{proof} According to Lemma~\ref{Mbounded}, the $\mathbf{A}$-vague adherence of
$(\mu_s)_{s\in S}$ is nonempty and contained in~$\mathfrak
M(\mathbf{A},\leqslant\!\mathbf{a},\mathbf{g})$. To establish the
lemma, it is enough to show that every  its element~$\mu$ is of
finite energy.

Observe that, by~(\ref{Ren}), the net of scalar measures
$(R\mu_s)_{s\in S}\subset\mathcal E$ is strongly bounded. We proceed
by proving that so are $(R\mu_s^+)_{s\in S}$ and $(R\mu_s^-)_{s\in
S}$, i.\,e.,
\begin{equation}
\sup_{s\in S}\,\|R\mu_s^\pm\|^2<\infty.  \label{7.1}
\end{equation}
Of course, this needs to be verified only when $I^-\ne\varnothing$;
then, according to the standing assumptions, both~(\ref{s})
and~(\ref{bou}) hold. Since $\langle g_i,\mu^i\rangle\leqslant a_i$,
we get
\begin{equation}
\sup_{s\in S}\,\mu_s^i(\mathbf X)\leqslant
a_ig_{i,\inf}^{-1}\quad\mbox{for all \ }i\in I. \label{7.3}
\end{equation}
Consequently, by (\ref{s}), \[\sup_{s\in S}\,R\mu_s^\pm(\mathbf
X)\leqslant \sum_{i\in I}\,a_ig_{i,\inf}^{-1}<\infty.\] Because
of~(\ref{bou}), this implies that $\kappa(R\mu^+_s,R\mu^-_s)$
remains bounded from above on~$S$; hence, so do $\|R\mu^+_s\|^2$ and
$\|R\mu^-_s\|^2$.

If $(\mu_d)_{d\in D}$ is a subnet of $(\mu_s)_{s\in S}$ that
converges $\mathbf{A}$-vaguely to~$\mu$, then, by
Lemma~\ref{lem:vague'}, $(R\mu^+_d)_{d\in D}$ and $(R\mu^-_d)_{d\in
D}$ converge vaguely to~$R\mu^+$ and~$R\mu^-$, respectively.
Therefore, applying Lemma~\ref{lemma:lower} with $\mathrm Y=\mathrm
X\times\mathrm X$ and $\psi=\kappa$, we conclude from~(\ref{7.1})
that $R\mu^+$ and $R\mu^-$ are both of finite energy. Because
of~(\ref{Ren}), this yields $\kappa(\mu,\mu)<\infty$, as was to be
proved.
\end{proof}

\begin{corollary}\label{cor:aux1} If a net $(\mu_s)_{s\in
S}\subset\mathcal E(\mathbf{A},\leqslant\!\mathbf{a},\mathbf{g})$ is
strongly bounded, then for every $i\in I$,
\begin{equation} \sup_{s\in S}\,\|\mu_s^i\|^2<\infty.\label{7.1i}
\end{equation}
\end{corollary}

\begin{proof} It is clear from~(\ref{7.1}) that the required relation
will be established once we prove
\begin{equation}\sum_{i,j\in
I^\pm}\,\kappa(\mu_s^i,\mu_s^j)\geqslant
C>-\infty,\label{bdbl}\end{equation} where $C$ is independent
of~$s$. Since (\ref{bdbl}) is obvious when $\kappa\geqslant0$, we
assume $\mathbf X$ to be compact. Then $\kappa$, being lower
semicontinuous, is bounded from below on~$\mathbf X$ (say by~$-c$,
where $c>0$), while $\mathbf{A}$ is finite. Furthermore, then
$g_{i,\inf}>0$; therefore, (\ref{7.3}) holds true. This implies that
$\kappa(\mu_s^i,\mu_s^j)\geqslant-a_ia_j\,g_{i,\inf}^{-1}\,g_{j,\inf}^{-1}\,c$
for all $i,\,j\in I$, and (\ref{bdbl}) follows.
\end{proof}

\subsection{Strong completeness of $\mathcal E(\mathbf{A},\leqslant\!\mathbf{a},\mathbf{g})$}\label{sec:strong}

\begin{theorem}\label{th:strong} The
semimetric space $\mathcal
E(\mathbf{A},\leqslant\!\mathbf{a},\mathbf{g})$ is complete. In more
detail, if $(\mu_s)_{s\in S}$ is a strong Cauchy net in $\mathcal
E(\mathbf{A},\leqslant\!\mathbf{a},\mathbf{g})$ and $\mu$~is its
$\mathbf{A}$-vague cluster point (such a~$\mu$ exists), then
$\mu\in\mathcal E(\mathbf{A},\leqslant\!\mathbf{a},\mathbf{g})$ and
\begin{equation}
\lim_{s\in S}\,\|\mu_s-\mu\|^2=0.\label{str}
\end{equation}
Assume, in addition, that the kernel $\kappa$ is strictly positive
definite and all $A_i$, $i\in I$, are mutually disjoint. If moreover
$(\mu_s)_{s\in S}\subset\mathcal
E(\mathbf{A},\leqslant\!\mathbf{a},\mathbf{g})$ converges strongly
to $\mu_0\in\mathcal E(\mathbf{A})$, then actually $\mu_0\in\mathcal
E(\mathbf{A},\leqslant\!\mathbf{a},\mathbf{g})$ and $\mu_s\to\mu_0$
$\mathbf{A}$-vaguely.\end{theorem}

\begin{proof}Fix a strong Cauchy net  $(\mu_s)_{s\in S}\subset\mathcal
E(\mathbf{A},\leqslant\!\mathbf{a},\mathbf{g})$. Since such a net
converges strongly to every its strong cluster point, $(\mu_s)_{s\in
S}$ can certainly be assumed to be strongly bounded. Then, by
Lemmas~\ref{Mbounded} and~\ref{lem:aux2}, there exists an
$\mathbf{A}$-vague cluster point~$\mu$ of~$(\mu_s)_{s\in S}$ and
\begin{equation}
\mu\in\mathcal E(\mathbf{A},\leqslant\!\mathbf{a},\mathbf{g}).
\label{leqslant1}
\end{equation}
We next proceed by verifying (\ref{str}). Of course, there is no
loss of generality in assuming $(\mu_s)_{s\in S}$ to converge
$\mathbf{A}$-vaguely to~$\mu$. Then, by Lemma~\ref{lem:vague'},
$(R\mu^+_s)_{s\in S}$ and $(R\mu^-_s)_{s\in S}$ converge vaguely
to~$R\mu^+$ and~$R\mu^-$, respectively. Since, by~(\ref{7.1}), these
nets are strongly bounded in~$\mathcal E^+$, the property~(C$_2$)
(see~Sect.~\ref{sec:2}) shows that they approach~$R\mu^+$
and~$R\mu^-$, respectively, in the weak topology as well, and so
$R\mu_s\to R\mu$ weakly. This gives, by~(\ref{seminorm}),
$$
\|\mu_s-\mu\|^2=\|R\mu_s-R\mu\|^2=\lim_{l\in
S}\,\kappa(R\mu_s-R\mu,R\mu_s-R\mu_l),
$$
and hence, by the Cauchy-Schwarz inequality,
$$\|\mu_s-\mu\|^2\leqslant
\|\mu_s-\mu\|\,\liminf_{l\in S}\,\|\mu_s-\mu_l\|,
$$
which proves (\ref{str}) as required, because $\|\mu_s-\mu_l\|$
becomes arbitrarily small when $s,\,l\in S$ are large enough.

Suppose now that $\kappa$ is strictly positive definite, while all
$A_i$, $i\in I$, are mutually disjoint, and let the net
$(\mu_s)_{s\in S}$ converge strongly to some $\mu_0\in\mathcal
E(\mathbf{A})$. Given an $\mathbf{A}$-vague limit point~$\mu$
of~$(\mu_s)_{s\in S}$, we conclude from~(\ref{str}) that
$\|\mu_0-\mu\|=0$, hence $R\mu_0=R\mu$ since $\kappa$ is strictly
positive definite, and finally $\mu_0=\mu$ because $A_i$, $i\in I$,
are mutually disjoint. In view of~(\ref{leqslant1}), this means that
$\mu_0\in\mathcal E(\mathbf{A},\leqslant\!\mathbf{a},\mathbf{g})$,
which is a part of the desired conclusion. Moreover, $\mu_0$ has
thus been shown to be identical to any $\mathbf{A}$-vague cluster
point of~$(\mu_s)_{s\in S}$. Since the $\mathbf{A}$-vague topology
is Hausdorff, this implies that $\mu_0$ is actually the
$\mathbf{A}$-vague limit of~$(\mu_s)_{s\in S}$ (cf.~\cite[Chap.~I,
\S~9, n$^\circ$\,1,
cor.]{B1}), which completes the proof.
\end{proof}

\begin{remark} In view of the fact that
the semimetric space $\mathcal
E(\mathbf{A},\leqslant\!\mathbf{a},\mathbf{g})$ is isometric to its
$R$-image, Theorem~\ref{th:strong} has thus singled out a {\it
strongly complete\/} topological subspace of the pre-Hilbert
space~$\mathcal E$, whose elements are {\it signed\/} measures. This
is of independent interest since, according to a well-known
counterexample by H.~Cartan~\cite{Car}, all the space~$\mathcal E$
is strongly incomplete even for the Newtonian kernel $|x-y|^{2-n}$
in~$\mathbb R^n$, $n\geqslant3$.\end{remark}

\begin{remark} Assume $\kappa$ is strictly positive
definite (hence, perfect). If moreover $I^-=\varnothing$, then
Theorem~\ref{th:strong} remains valid for $\mathcal E(\mathbf{A})$
in place of $\mathcal
E(\mathbf{A},\leqslant\!\mathbf{a},\mathbf{g})$
(cf.~Theorem~\ref{th:1}). A question still unanswered is whether
this is the case if $I^+$ and $I^-$ are both nonempty. We can
however show that this is really so for the Riesz kernels
$|x-y|^{\alpha-n}$, $0<\alpha<n$, in~$\mathbb R^n$, $n\geqslant2$
(cf.~\cite[Th.~1]{Z2}). The proof utilizes Deny's theorem~\cite{D1}
stating that, for the Riesz kernels, $\mathcal E$~can be completed
with making use of distributions of finite energy.\end{remark}

\section{Extremal measures in the Gauss variational problem}\label{sec:extremal}
To apply Theorem~\ref{th:strong} to the Gauss variational problem,
we next proceed by introducing the concept of extremal measure
defined as a strong and, simultaneously, the $\mathbf{A}$-vague
limit of a minimizing net. See below for strict definitions and
related auxiliary results.

Except for Corollary~\ref{IorII}, in addition to the standing
assumptions we suppose that
\begin{equation}\label{minfty}G_{\mathbf{f}}(\mathbf{A},\mathbf{a},\mathbf{g})>-\infty.\end{equation}

\subsection{Extremal measures: existence, uniqueness, and $\mathbf{A}$-vague compactness}

\begin{definition}We call a net $(\mu_s)_{s\in S}$ {\it minimizing\/} if
$(\mu_s)_{s\in S}\subset\mathcal
E^0_{\mathbf{f}}(\mathbf{A},\mathbf{a},\mathbf{g})$ and
\begin{equation}\lim_{s\in S}\,G_{\mathbf{f}}(\mu_s)=G_{\mathbf{f}}(\mathbf{A},\mathbf{a},\mathbf{g}).
\label{min}\end{equation}\end{definition}

Let $\mathbb M_{\mathbf{f}}(\mathbf{A},\mathbf{a},\mathbf{g})$
consist of all minimizing nets; note that it is nonempty, which is
clear from~(\ref{nonzero1}) and~Corollary~\ref{cor:compact}. We
denote by $\mathcal
M_{\mathbf{f}}(\mathbf{A},\mathbf{a},\mathbf{g})$ the union of the
$\mathbf{A}$-vague cluster sets of $(\mu_s)_{s\in S}$, where
$(\mu_s)_{s\in S}$ ranges over $\mathbb
M_{\mathbf{f}}(\mathbf{A},\mathbf{a},\mathbf{g})$.

\begin{definition}\label{def:extr} We call $\gamma\in\mathcal E(\mathbf{A})$ {\it extremal\/}
if there exists $(\mu_s)_{s\in S}\in\mathbb
M_{\mathbf{f}}(\mathbf{A},\mathbf{a},\mathbf{g})$ that converges
to~$\gamma$ both strongly and $\mathbf{A}$-vaguely; such a net
$(\mu_s)_{s\in S}$ is said to {\it generate\/}~$\gamma$. The class
of all extremal measures will be denoted by $\mathfrak
E_{\mathbf{f}}(\mathbf{A},\mathbf{a},\mathbf{g})$.
\end{definition}

\begin{lemma}\label{lemma:WM} The following assertions hold true:
\begin{itemize}
\item[\rm(i)] From every minimizing net one can select a subnet generating an extremal measure;
hence, $\mathfrak E_{\mathbf{f}}(\mathbf{A},\mathbf{a},\mathbf{g})$
is nonempty. Furthermore,
\begin{equation}\label{ME}\mathfrak E_{\mathbf{f}}(\mathbf{A},\mathbf{a},\mathbf{g})\subset\mathcal
E(\mathbf{A},\leqslant\!\mathbf{a},\mathbf{g})\end{equation} and
\begin{equation}\mathfrak E_{\mathbf{f}}(\mathbf{A},\mathbf{a},\mathbf{g})=\mathcal M_{\mathbf{f}}(\mathbf{A},\mathbf{a},\mathbf{g}).
\label{WM}\end{equation}
\item[\rm(ii)] Every minimizing
net converges strongly to every extremal measure; hence, $\mathfrak
E_{\mathbf{f}}(\mathbf{A},\mathbf{a},\mathbf{g})$ is contained in an
equivalence class in~$\mathcal E(\mathbf{A})$.\smallskip
\item[\rm(iii)] The class $\mathfrak E_{\mathbf{f}}(\mathbf{A},\mathbf{a},\mathbf{g})$ is $\mathbf{A}$-vaguely compact.
\end{itemize}
\end{lemma}

\begin{proof}Fix $(\mu_s)_{s\in S}$ and $(\nu_t)_{t\in T}$ in $\mathbb
M_{\mathbf{f}}(\mathbf{A},\mathbf{a},\mathbf{g})$. Then
\begin{equation}
\lim_{(s,\,t)\in S\times T}\,\|\mu_s-\nu_t\|^2=0, \label{fund}
\end{equation}
where $S\times T$ denotes the directed product of the directed
sets~$S$ and~$T$ (see, e.\,g.,~\cite[Chap.~2,~\S~3]{K}). Indeed,
since $\mathcal E_{\mathbf{f}}^0(\mathbf{A},\mathbf{a},\mathbf{g})$
is convex, in the same manner as in the proof
of~Lemma~\ref{lemma:unique} we get
\[0\leqslant\|R\mu_s-R\nu_t\|^2\leqslant-4G_{\mathbf{f}}(\mathbf{A},\mathbf{a},\mathbf{g})+
2G_{\mathbf{f}}(\mu_s)+2G_{\mathbf{f}}(\nu_t),\] which yields
(\ref{fund}) when combined with (\ref{min}).

Relation~(\ref{fund}) implies that $(\mu_s)_{s\in S}$ is strongly
fundamental. Therefore, by Theorem~\ref{th:strong}, there is an
$\mathbf{A}$-vague cluster point~$\mu_0$ of~$(\mu_s)_{s\in S}$,
$\mu_0\in\mathcal E(\mathbf{A},\leqslant\!\mathbf{a},\mathbf{g})$,
and $\mu_s\to\mu_0$ strongly. This means that $\mu_0$ is an extremal
measure and, hence, $\mathcal
M_{\mathbf{f}}(\mathbf{A},\mathbf{a},\mathbf{g})\subset\mathfrak
E_{\mathbf{f}}(\mathbf{A},\mathbf{a},\mathbf{g})$. Since the inverse
inclusion is obvious, relations~(\ref{ME}) and~(\ref{WM}) follow.

To verify (ii), fix $(\mu_s)_{s\in S}\in\mathbb
M_{\mathbf{f}}(\mathbf{A},\mathbf{a},\mathbf{g})$ and
$\gamma\in\mathfrak
E_{\mathbf{f}}(\mathbf{A},\mathbf{a},\mathbf{g})$. Then, by
Definition~\ref{def:extr}, one can choose a net in~$\mathbb
M_{\mathbf{f}}(\mathbf{A},\mathbf{a},\mathbf{g})$, say
$(\nu_t)_{t\in T}$, that converges to~$\gamma$ strongly. Repeated
application of~(\ref{fund}) shows that also $(\mu_s)_{s\in S}$
converges to~$\gamma$ strongly, as claimed.

To establish (iii), it is enough to prove that $\mathcal
M_{\mathbf{f}}(\mathbf{A},\mathbf{a},\mathbf{g})$ is
$\mathbf{A}$-vaguely compact. Fix $(\gamma_s)_{s\in
S}\subset\mathcal M_{\mathbf{f}}(\mathbf{A},\mathbf{a},\mathbf{g})$.
It follows from~(\ref{ME}) and~Lemma~\ref{Mbounded} that there
exists an $\mathbf{A}$-vague cluster point~$\gamma_0$ of
$(\gamma_s)_{s\in S}$; let $(\gamma_t)_{t\in T}$ be a subnet
of~$(\gamma_s)_{s\in S}$ that converges $\mathbf{A}$-vaguely
to~$\gamma_0$. Then for every $t\in T$ one can choose
$(\mu_{s_t})_{s_t\in S_t}\in\mathbb
M_{\mathbf{f}}(\mathbf{A},\mathbf{a},\mathbf{g})$ converging
$\mathbf{A}$-vaguely to~$\gamma_t$. Consider the Cartesian product
$\prod\,\{S_t: t\in T\}$~--- that is, the collection of all
functions~$\beta$ on~$T$ with $\beta(t)\in S_t$, and let~$D$ denote
the directed product $T\times\prod\,\{S_t: t\in T\}$. Given
$(t,\beta)\in D$, write $\mu_{(t,\,\beta)}:=\mu_{\beta(t)}$. Then
the theorem on iterated limits from \cite[Chap.~2, \S~4]{K} yields
that the net $(\mu_{(t,\beta)})_{(t,\beta)\in D}$ belongs to
$\mathbb M_{\mathbf{f}}(\mathbf{A},\mathbf{a},\mathbf{g})$ and
converges $\mathbf{A}$-vaguely to~$\gamma_0$. Thus,
$\gamma_0\in\mathcal
M_{\mathbf{f}}(\mathbf{A},\mathbf{a},\mathbf{g})$ as was to be
proved.
\end{proof}

\begin{corollary}\label{cor:WS}
Every equilibrium measure $\lambda$ (if exists) is extremal, i.\,e.,
\begin{equation}\label{WS}\mathcal G_{\mathbf{f}}(\mathbf{A},\mathbf{a},\mathbf{g})\subset\mathfrak
E_{\mathbf{f}}(\mathbf{A},\mathbf{a},\mathbf{g}).\end{equation} If
$(\mu_s)_{s\in S}\in\mathbb
M_{\mathbf{f}}(\mathbf{A},\mathbf{a},\mathbf{g})$ is arbitrarily
given, then $\mu_s\to\lambda$ strongly and, moreover,
\begin{equation}\label{flambda}\lim_{s\in S}\,\langle
f,\mu_s\rangle=\langle f,\lambda\rangle.\end{equation}
\end{corollary}

\begin{proof}
For every $\mathbf{K}\in\{\mathbf{K}\}_{\mathbf{A}}$ large enough
consider
$\hat{\lambda}_{\mathbf{K}}:=\bigl(\hat{\lambda}_{\mathbf{K}}^i\bigr)_{i\in
I}$, where $\hat{\lambda}_{\mathbf{K}}^i$ is given by~(\ref{hatmu})
with $\mu=\lambda$. Then
$(\hat{\lambda}_{\mathbf{K}})_{\mathbf{K}\in\{\mathbf{K}\}_{\mathbf{A}}}$
belongs to $\mathbb
M_{\mathbf{f}}(\mathbf{A},\mathbf{a},\mathbf{g})$, which is clear
from~(\ref{4w}) with~$\mu$ replaced by~$\lambda$. On the other hand,
this net converges $\mathbf{A}$-vaguely to~$\lambda$; hence,
$\lambda\in\mathcal
M_{\mathbf{f}}(\mathbf{A},\mathbf{a},\mathbf{g})$. Therefore, in
accordance with~(\ref{WM}), $\lambda$ has to be extremal.

Fix $(\mu_s)_{s\in S}\in\mathbb
M_{\mathbf{f}}(\mathbf{A},\mathbf{a},\mathbf{g})$; then
$\mu_s\to\lambda$ strongly, which is a consequence of~(\ref{WS}) and
Lemma~\ref{lemma:WM},~(ii). This implies that $\lim_{s\in
S}\,\|\mu_s\|^2=\|\lambda\|^2$. On the other hand, by~(\ref{min}),
\[\|\lambda\|^2+2\langle
f,\lambda\rangle=G_{\mathbf{f}}(\mathbf{A},\mathbf{a},\mathbf{g})=\lim_{s\in
S}\,\bigl[\|\mu_s\|^2+2\langle f,\mu_s\rangle\bigr].\] The last two
relations combined give~(\ref{flambda}), and the proof is
complete.
\end{proof}

\begin{corollary}\label{IorII}Assume that Case I or II occurs. Then
$G_{\mathbf{f}}(\mathbf{A},\mathbf{a},\mathbf{g})>-\infty$ and,
moreover,
\begin{equation}\label{conclII}G_{\mathbf{f}}(\gamma)=G_{\mathbf{f}}(\mathbf{A},\mathbf{a},\mathbf{g})\quad
\mbox{for all \ }\gamma\in\mathfrak
E_{\mathbf{f}}(\mathbf{A},\mathbf{a},\mathbf{g}).\end{equation}
\end{corollary}

\begin{proof}Suppose Case~II takes place; then $f_i=\alpha_i\kappa(\,\cdot\,,\sigma)$ for all $i\in I$, where
$\sigma\in\mathcal E$. Hence,
\begin{equation*}\label{II}
\langle\mathbf{f},\mu\rangle=\sum_{i\in
I}\,\alpha_i\int\kappa(x,\sigma)\,d\mu^i(x)=\kappa(\sigma,R\mu)\quad\mbox{for
all \ } \mu\in\mathcal E(\mathbf{A}),
\end{equation*}
the latter equality being a consequence of Lemma~\ref{integral}.
This implies
\begin{equation}\label{yus}
G_{\mathbf{f}}(\mu)=\|\mu\|^2+2\kappa(\sigma,R\mu)=\|R\mu+\sigma\|^2-\|\sigma\|^2.\end{equation}
Therefore
$G_{\mathbf{f}}(\mathbf{A},\mathbf{a},\mathbf{g})\geqslant-\|\sigma\|^2>-\infty$,
which enables us to use Lemma~\ref{lemma:WM}.

Applying (\ref{yus}) to $\mu_s$, $s\in S$, and $\gamma$, where
$(\mu_s)_{s\in S}\in\mathbb
M_{\mathbf{f}}(\mathbf{A},\mathbf{a},\mathbf{g})$ and
$\gamma\in\mathfrak
E_{\mathbf{f}}(\mathbf{A},\mathbf{a},\mathbf{g})$ are arbitrarily
given, in view of the fact that $\mu_s\to\gamma$ strongly we get
\[G_{\mathbf{f}}(\gamma)=\|R\gamma+\sigma\|^2-\|\sigma\|^2=\lim_{s\in S}\,
\bigl[\|R\mu_s+\sigma\|^2-\|\sigma\|^2\bigr]=\lim_{s\in
S}\,G_{\mathbf{f}}(\mu_s).\] Substituting (\ref{min}) into the
preceding relation yields (\ref{conclII}).

Since, by~(\ref{Rpot}), Case~I can be reduced to Case~II with
$\sigma=R\nu$, the proof is complete.
\end{proof}

\subsection{Extremal measures: $g_i$-masses of the $i$-components}

\begin{lemma}\label{lemma:exist}
Fix $i\in I$ and assume that either $g_{i,\sup}<\infty$ or
{\rm(\ref{growth})} holds for some $r_i\in(1,\infty)$ and
$\omega_i\in\mathcal E$. If moreover $A_i$ either is compact or has
finite interior capacity, then
\begin{equation}
\langle g_i,\gamma^i\rangle=a_i\quad\mbox{for all \
}\gamma\in\mathfrak
E_{\mathbf{f}}(\mathbf{A},\mathbf{a},\mathbf{g}). \label{24}
\end{equation}\end{lemma}

\begin{proof}Fix $\gamma\in\mathfrak E_{\mathbf{f}}(\mathbf{A},\mathbf{a},\mathbf{g})$ and choose $(\mu_s)_{s\in
S}\in\mathbb M_{\mathbf{f}}(\mathbf{A},\mathbf{a},\mathbf{g})$
generating~$\gamma$. Taking a subnet if necessary, one can assume
$(\mu_s)_{s\in S}$ to be strongly bounded. Then, by~(\ref{7.1i}), so
is~$(\mu^i_s)_{s\in S}$.

Of course, (\ref{24}) needs to be proved only if the set~$A_i$ is
noncompact; then its capacity has to be finite. Hence,
by~\cite[Th.~4.1]{F1}, for every $E\subset A_i$ there exists a
measure $\theta_E\in\mathcal E^+(\,\overline{E}\,)$, called an
interior equilibrium measure associated with~$E$, which possesses
the properties
\begin{equation}
\theta_E(\mathrm X)=\|\theta_E\|^2=C(E), \label{5}
\end{equation}
\begin{equation}
\kappa(x,\theta_E)\geqslant1\quad\mbox{n.\,e. in \ } E. \label{6}
\end{equation}

Also observe that there is no loss of generality in assuming $g_i$
to satisfy~(\ref{growth}) with some $r_i\in(1,\infty)$ and
$\omega_i\in\mathcal E$. Indeed, otherwise $g_i$ has to be bounded
from above (say by $M$), which combined with~(\ref{6}) again
gives~(\ref{growth}) for $\omega_i:=M^{r_i}\,\theta_{A_i}$,
$r_i\in(1,\infty)$ being arbitrary.

To establish (\ref{24}), we treat $A_i$ as a locally compact space
with the topology induced from~$\mathrm X$. Given a set $E\subset
A_i$, let $\chi_E$ denote its characteristic function and let
$E^c:=A_i\setminus E$. Further, let $\{K_i\}$ be the increasing
family of all compact subsets~$K_i$ of~$A_i$. Since $g_i\chi_{K_i}$
is upper semicontinuous on~$A_i$ while $(\mu_s^i)_{s\in S}$
converges to~$\gamma^i$ vaguely, for every $K_i\in\{K_i\}$
\[
\langle g_i\chi_{K_i},\gamma^i\rangle\geqslant\limsup_{s\in
S}\,\langle g_i\chi_{K_i},\mu_s^i\rangle\] according to
Lemma~\ref{lemma:lower}. On the other hand, application of
Lemma~1.2.2 from~\cite{F1} yields
\[
\langle g_i,\gamma^i\rangle=\lim_{K_i\in\{K_i\}}\,\langle
g_i\chi_{K_i},\gamma^i\rangle.\] Combining the last two relations,
we obtain
\[
a_i\geqslant\langle
g_i,\gamma^i\rangle\geqslant\limsup_{(s,\,K_i)\in
S\times\{K_i\}}\,\langle g_i\chi_{K_i},\mu_s^i\rangle=
a_i-\liminf_{(s,\,K_i)\in S\times\{K_i\}}\,\langle
g_i\chi_{K^c_i},\mu_s^i\rangle,\] $S\times\{K_i\}$ being the
directed product of the directed sets~$S$ and~$\{K_i\}$. Hence, if
we prove
\begin{equation}
\liminf_{(s,\,K_i)\in S\times\{K_i\}}\,\langle
g_i\chi_{K^c_i},\mu_s^i\rangle=0, \label{25}
\end{equation}
the desired relation~(\ref{24}) follows.

Consider an interior equilibrium measure $\theta_{K^c_i}$, where
$K_i\in\{K_i\}$ is given. Then application of~Lemma~4.1.1 and
Theorem~4.1 from~\cite{F1} shows that
\[
\|\theta_{K^c_i}-\theta_{\tilde{K}^c_i}\|^2\leqslant
\|\theta_{K^c_i}\|^2-\|\theta_{\tilde{K}^c_i}\|^2\quad\mbox{provided
\ }K_i\subset\tilde{K}_i.\] Furthermore, it is clear from~(\ref{5})
that the net $\|\theta_{K^c_i}\|$, $K_i\in\{K_i\}$, is bounded and
nonincreasing, and hence fundamental in~$\mathbb R$. The preceding
inequality thus yields that the net
$(\theta_{K^c_i})_{K_i\in\{K_i\}}$ is strongly fundamental
in~$\mathcal E$. Since, clearly, it converges vaguely to zero, the
property~(C$_1$) (see.~Sec.~\ref{sec:2}) implies immediately that
zero is also one of its strong limits and, hence,
\begin{equation}
\lim_{K_i\in\{K_i\}}\,\|\theta_{K^c_i}\|=0. \label{27}
\end{equation}

Write $q_i:=r_i(r_i-1)^{-1}$, where $r_i\in(1,\infty)$ is a number
involved in condition~(\ref{growth}). Combining (\ref{growth}) with
(\ref{6}) shows that the inequality
\[
g_i(x)\,\chi_{K^c_i}(x)\leqslant\kappa(x,\omega_i)^{1/r_i}\,
\kappa(x,\theta_{K^c_i})^{1/q_i}\]  subsists n.\,e.~in~$A_i$, and
hence $\mu_s^i$-a.\,e.~in~$\mathrm X$ by virtue
of~\cite[Lemma~2.3.1]{F1} and the fact that $\mu_s^i$ is a measure
of finite energy, compactly supported in~$A_i$. Having integrated
this relation with respect to~$\mu_s^i$, we then apply the H\"older
and, subsequently, the Cauchy-Schwarz inequalities to the integrals
on the right. This gives
\[\langle
g_i\chi_{K^c_i},\mu_s^i\rangle\leqslant\Bigl[\int\kappa(x,\omega_i)\,d\mu_s^i(x)\Bigr]^{1/r_i}\,
\Bigl[\int\kappa(x,\theta_{K^c_i})\,d\mu_s^i(x)\Bigr]^{1/q_i}\leqslant
\|\omega_i\|^{1/r_i}\,\|\theta_{K^c_i}\|^{1/q_i}\,\|\mu_s^i\|.\]
Taking limits here along $S\times\{K\}$ and using (\ref{7.1i}) and
(\ref{27}), we obtain~(\ref{25}) as desired.
\end{proof}

\section{Proof of Theorem~\ref{exist}}\label{sec:proof.th.str}

We begin by verifying relation~(\ref{minfty}). This needs to be done
only in Case~III, because in the remaining Cases~I and~II it has
already been established by Corollary~\ref{IorII}. In view of the
positive definiteness of the kernel, it suffices to show that
\begin{equation}\label{minus}\langle\mathbf{f},\mu\rangle\geqslant-M_0>-\infty\quad\mbox{for all
\ } \mu\in\mathcal
E(\mathbf{A},\mathbf{a},\mathbf{g}).\end{equation} Assume $\mathrm
X$ to be compact, since otherwise $f_i\geqslant0$ for all $i\in I$
and (\ref{minus}) is obvious. Then $\mathbf{A}$ is finite and, for
every $i\in I$, $g_{i,\inf}>0$ while $f_i$, being lower
semicontinuous, is bounded from below. This implies~(\ref{minus})
when combined with the inequalities $\mu^i(\mathrm X)\leqslant
a_i\,g_{i,\inf}^{-1}<\infty$.

Due to (\ref{minfty}), we are able to use the results from
Sect.~\ref{sec:extremal}. Fix an extremal measure~$\gamma$~--- it
exists according to~Lemma~\ref{lemma:WM}, and choose a net
$(\mu_s)_{s\in S}\in\mathbb
M_{\mathbf{f}}(\mathbf{A},\mathbf{a},\mathbf{g})$ that converges
to~$\gamma$ both strongly and $\mathbf{A}$-vaguely. We are going to
prove that $\gamma$ is an equilibrium measure.

Observe that, by~Lemma~\ref{lemma:exist}, $\gamma\in\mathcal
E(\mathbf{A},\mathbf{a},\mathbf{g})$. Hence, the desired inclusion
$\gamma\in\mathcal G_{\mathbf{f}}(\mathbf{A},\mathbf{a},\mathbf{g})$
will have been established once we show that
$\langle\mathbf{f},\gamma\rangle>-\infty$ and
\begin{equation}\label{minimizer}G_{\mathbf{f}}(\gamma)\leqslant
G_{\mathbf{f}}(\mathbf{A},\mathbf{a},\mathbf{g}).\end{equation} To
this end, one can again assume Case~III to occur, for otherwise this
has already been obtained by Corollary~\ref{IorII}. Then
$\langle\mathbf{f},\gamma\rangle>-\infty$ by~(\ref{minus}) for
$\gamma$ instead of~$\mu$. Furthermore, from the strong and the
$\mathbf{A}$-vague convergence of~$(\mu_s)_{s\in S}$ to~$\gamma$ we
respectively get
\[G_{\mathbf{f}}(\mathbf{A},\mathbf{a},\mathbf{g})=\lim_{s\in
S}\,\bigl[\|\mu_s\|^2+2\langle
f,\mu_s\rangle\bigr]=\|\gamma\|^2+2\lim_{s\in S}\,\langle
f,\mu_s\rangle\] and
\begin{equation*}\label{hryu}\sum_{i\in I}\,\langle f_i,\gamma^i\rangle\leqslant
\sum_{i\in I}\,\liminf_{s\in S}\,\langle
f_i,\mu_s^i\rangle\leqslant\lim_{s\in S}\,\sum_{i\in I}\,\langle
f_i,\mu_s^i\rangle.\end{equation*}  The last two relations combined
give~(\ref{minimizer}).

What has thus been proved means that the Gauss variational problem
is solvable; actually, $\mathfrak
E_{\mathbf{f}}(\mathbf{A},\mathbf{a},\mathbf{g})\subset\mathcal
G_{\mathbf{f}}(\mathbf{A},\mathbf{a},\mathbf{g})$. Together
with~(\ref{WM}) and~(\ref{WS}), this yields
\begin{equation}\label{SWM}\mathcal G_{\mathbf{f}}(\mathbf{A},\mathbf{a},\mathbf{g})=
\mathfrak E_{\mathbf{f}}(\mathbf{A},\mathbf{a},\mathbf{g})=\mathcal
M_{\mathbf{f}}(\mathbf{A},\mathbf{a},\mathbf{g}).\end{equation}
Therefore Lemma~\ref{lemma:WM}, (iii), implies that $\mathcal
G_{\mathbf{f}}(\mathbf{A},\mathbf{a},\mathbf{g})$ is
$\mathbf{A}$-vaguely
compact.\hfill$\square$

\section{Proof of Theorem~\ref{cor:cont}}\label{sec:proof.th.cont}

Fix $\lambda_{\mathbf{K}}\in\mathcal
G_{\mathbf{f}}(\mathbf{K},\mathbf{a},\mathbf{g})$, where
$\mathbf{K}\in\{\mathbf{K}\}_{\mathbf{A}}$, and
$\lambda_{\mathbf{A}}\in\mathcal
G_{\mathbf{f}}(\mathbf{A},\mathbf{a},\mathbf{g})$
--- the existence of such equilibrium measures has been ensured
by~Theorem~\ref{exist} and~Corollary~\ref{cor:exist}. According
to~Lemma~\ref{lemma.cont},
\begin{equation}\label{mincom}(\lambda_{\mathbf{K}})_{\mathbf{K}\in\{\mathbf{K}\}_{\mathbf{A}}}\in\mathbb
M_{\mathbf{f}}(\mathbf{A},\mathbf{a},\mathbf{g}).\end{equation}
Therefore, by~(\ref{SWM}), every $\mathbf{A}$-vague cluster point of
$(\lambda_{\mathbf{K}})_{\mathbf{K}\in\{\mathbf{K}\}_{\mathbf{A}}}$
belongs to~$\mathcal
G_{\mathbf{f}}(\mathbf{A},\mathbf{a},\mathbf{g})$, which is a part
of the desired conclusion. Furthermore, the claimed
relations~(\ref{cor:cont1}) and~(\ref{cor:cont2}) are obtained
directly from~(\ref{mincom}) and~Corollary~\ref{cor:WS}. What is
thus left is to establish~(\ref{cor:cont3}).

Consider an arbitrary cluster point $d_i$ of $\langle
W^i_{\lambda_{\mathbf{K}}},\lambda_{\mathbf{K}}^i\rangle$, where
$\mathbf{K}$ ranges over $\{\mathbf{K}\}_{\mathbf{A}}$. Then
application of Lemma\ref{lemma:WM},~(i), implies that there exists a
subnet $(\lambda_s)_{s\in S}$ of
$(\lambda_{\mathbf{K}})_{\mathbf{K}\in\{\mathbf{K}\}_{\mathbf{A}}}$,
strongly and $\mathbf{A}$-vaguely convergent (say to~$\lambda$) and
such that
\begin{equation}\label{D}d_i=\lim_{s\in S}\,\langle
W^i_{\lambda_s},\lambda_s^i\rangle.\end{equation} Also observe that,
by (\ref{SWM}) and (\ref{mincom}), $\lambda\in\mathcal
G_{\mathbf{f}}(\mathbf{A},\mathbf{a},\mathbf{g})$; hence,
by~Corollary~\ref{cor:unique}, \begin{equation}\label{DD}\langle
W^i_{\lambda},\lambda^i\rangle=\langle
W^i_{\lambda_{\mathbf{A}}},\lambda_{\mathbf{A}}^i\rangle.\end{equation}
We proceed by showing that, for every $i\in I$,
\begin{align}\label{stepI1}\langle\kappa^i_\lambda,\lambda^i\rangle&=
\lim_{s\in S}\,\langle\kappa^i_{\lambda_s},\lambda_s^i\rangle,\\
\label{stepI2}\langle f_i,\lambda^i\rangle&=\lim_{s\in S}\,\langle
f_i,\lambda_s^i\rangle.\end{align}

Without loss of generality, $(\lambda_s)_{s\in S}$ can certainly be
assumed to be strongly bounded. Then, by Corollary~\ref{cor:aux1},
so is $(\lambda^i_s)_{s\in S}$. Since, moreover,
$\lambda^i_s\to\lambda^i$ vaguely, the property~(C$_2$) implies that
$\lambda^i_s$ approaches~$\lambda^i$ also weakly. Hence, for every
$\varepsilon>0$,
$\bigl|\kappa(\lambda^i-\lambda^i_s,R\lambda)\bigr|<\varepsilon$
whenever $s\in S$ is large enough. Furthermore, by the
Cauchy-Schwarz inequality,
\[\bigl|\kappa(\lambda^i_s,R\lambda)-\kappa(\lambda^i_s,R\lambda_s)\bigr|=
\bigl|\kappa(\lambda^i_s,R\lambda-R\lambda_s)\bigr|\leqslant
M_1\,\|\lambda-\lambda_s\|^2,\quad s\in S.\] Since
$\lambda_s\to\lambda$ strongly, the last two relations combined
yield
\[\kappa(\lambda^i,R\lambda)=\lim_{s\in S}\,\kappa(\lambda_s^i,R\lambda_s),\]
which in view of~(\ref{Rpot}) is equivalent to~(\ref{stepI1}).

To establish~(\ref{stepI2}), we can restrict ourselves to Case~III,
for otherwise it is obtained directly from the weak convergence
of~$(\lambda^i_s)_{s\in S}$ to~$\lambda^i$. Then it follows from the
$\mathbf{A}$-vague convergence of~$(\lambda_s)_{s\in S}$
to~$\lambda$ that
\begin{equation}\label{lastbut} \langle f_i,\lambda^i\rangle\leqslant\liminf_{s\in S}\,\langle
f_i,\lambda_s^i\rangle\quad\mbox{for all \ } i\in I\end{equation}
and therefore, by (\ref{flambda}),
\[\langle\mathbf{f},\lambda\rangle=\sum_{i\in I}\,\langle f_i,\lambda^i\rangle\leqslant
\sum_{i\in I}\,\liminf_{s\in S}\,\langle
f_i,\lambda_s^i\rangle\leqslant\lim_{s\in S}\,\langle
\mathbf{f},\lambda_s\rangle=\langle\mathbf{f},\lambda\rangle.\]
Comparing the last two relations yields that an equality
in~(\ref{lastbut}) actually has to hold. Since the same arguments
can be applied to any subnet of~$(\lambda_s)_{s\in S}$,
(\ref{stepI2}) follows.

Combining~\mbox{(\ref{D})--(\ref{stepI2})} shows that $d_i=\langle
W^i_{\lambda_{\mathbf{A}}},\lambda_{\mathbf{A}}^i\rangle$. Since
this has been established for any cluster point~$d_i$ of $\langle
W^i_{\lambda_{\mathbf{K}}},\lambda_{\mathbf{K}}^i\rangle$, where
$\mathbf{K}\in\{\mathbf{K}\}_{\mathbf{A}}$, the claimed
relation~(\ref{cor:cont3}) is proved.\hfill$\square$

\section{Acknowledgments} The author wishes to thank Professors A.\,Yu.~Rashkovskii and
W.\,L.~Wendland for the careful reading of the manuscript and their
suggestions to improve it.

A part of this research was done during the author's visits to the
University Stuttgart and the Mathematisches Forschungsinstitut
Oberwolfach during April--May of 2009, and the author acknowledges
these institutions for the support and the excellent working
conditions.

\end{document}